\documentclass[
final,
nomarks
]{dmtcs-episciences}

\usepackage[utf8]{inputenc}
\usepackage{subfigure}
\usepackage{ amssymb }

\usepackage{graphicx}

\usepackage[round]{natbib}

\author[Samvel Kh. Darbinyan]{Samvel Kh. Darbinyan
 }
\title[A new sufficient condition 
 for a 2-strong digraph to be   Hamiltonian]{A new sufficient condition   for a 2-strong digraph to be   Hamiltonian}

\affiliation{ Institute for Informatics and Automation Problems of NAS RA,  Yerevan, Armenia}

\begin{document}
% This is only used if you are compiling for a volume before vol 25
% \publicationdetails{VOL}{2015}{ISS}{NUM}{SUBM}
% This is the new form of collecting the data, starting with vol 25

\publicationdata
{vol. 26:2}
{2024}
{7}
{10.46298/dmtcs.11560}
%{1998-10-14; 1998-10-14; 2002-07-19; 2014-02-05; 2015-09-09; 2022-12-25}
%{2022-12-3}
%{2023-12-3; None}
%{2023-1-1}
{2023-07-10; 2023-07-10; 2024-01-09}
{2024-03-25}

\maketitle

\begin{abstract}

 In this paper we prove the following new sufficient condition for a digraph to be Hamiltonian: 

{\it Let $D$ be a 2-strong  digraph of order $n\geq 9$. 
  If $n-1$ vertices of $D$ have degrees at least $n+k$ and the remaining
  vertex has degree at least $n-k-4$, where $k$ is a non-negative  integer, then $D$  is Hamiltonian}. 
 
 This  is  an extension of Ghouila-Houri's theorem for 2-strong digraphs and is a generalization of an early result of the author (DAN Arm. SSR (91(2):6-8, 1990). The obtained result is best possible in the sense that for $k=0$ there is a digraph of order $n=8$ (respectively, $n=9$) with the minimum degree $n-4=4$ (respectively, with the minimum degree $n-5=4$) whose $n-1$ vertices have degrees at least $n-1$, but it is not Hamiltonian. 

 We also give a new sufficient condition for a 3-strong digraph to be Hamiltonian-connected.
\end{abstract}

%DMTCS is an open access scientific is implemented by the
%\emph{episcience} platform, see \cite{berthaud:hal-01002815} for an
%overview of the strategy. It combines high scientific and editorial
%quality with an open access policy. It is priceless, neither authors
%nor readers pay money for the access. Access is granted by giving
%episcience an irrevocable license to publish the articles, the
%copyright remains with the authors. The platform itself is run by
%French government services that do their best to warrant continuous
%access and a high quality of service.

%This document describes the use of the \texttt{dmtcs-episcience.cls}
%document class. It has to be used \emph{for all DMTCS publications}.

\section{Introduction}

 In this paper, we consider finite digraphs without loops and multiple arcs. We shall assume that the reader is familiar with the standard  terminology on digraphs and refer the reader to \cite{[5]}. Every cycle and path is assumed simple and directed. A cycle (a path) in a digraph $D$ is called {\it Hamiltonian} ({\it Hamiltonian path}) if it  includes all the vertices of $D$. A digraph $D$ is {\it Hamiltonian}  if it contains a  Hamiltonian cycle. Hamiltonicity is one of the most central in graph theory, and it has been extensively studied by numerous researchers. The problem of deciding Hamiltonicity of a graph (digraph) is $NP$-complete, but there  
  are numerous sufficient conditions which ensure  the existence of a Hamiltonian cycle in a digraph (see \cite{[5], [6], [16], [18]}). Among them 
  are  the following classical sufficient conditions for a digraph to be Hamiltonian.\\
 
 \noindent\textbf{Theorem 1.1} (\cite{[20]}).  {\it Let $D$ be a  digraph of order $n\geq 2$. If for every vertex $x$ of $D$, 
 $d^+(x)\geq n/2$ and $d^-(x)\geq n/2$, then $D$ is Hamiltonian.} \\
 
 \noindent\textbf{Theorem 1.2} (\cite{[14]}). {\it Let $D$ be a strong digraph of order $n\geq 2$. If for every vertex $x$ of $D$, $d(x)\geq n$, then $D$ is Hamiltonian.}\\ 
 
  \noindent\textbf{Theorem 1.3} (\cite{[27]}).  {\it Let $D$ be  digraph of order $n\geq 2$. If  $d^+(x)+d^-(y)\geq n$ for all pairs of distinct vertices $x$ and $y$ of $D$ such that there is no arc from $x$, to $y$, then $D$ is Hamiltonian.}\\ 
  
 \noindent\textbf{Theorem 1.4} (\cite{[19]}).  Let $D$ be  a strong digraph of order $n\geq 2$. If $d(x)+d(y)\geq 2n-1$ for all pairs of  non-adjacent distinct  vertices  $x$ and $y$ of $D$, then $D$ is Hamiltonian.\\

It is known that all the lower bounds in the above theorems are tight.
Notice that for the strong digraphs Meyniel's theorem is a  generalization of Nash-Williams, Ghouila-Houri's and Woodall's theorems. A beautiful short proof the later can found in the paper by \cite{[7]}.

  \cite{[20]} suggested the problem of characterizing all the strong  digraphs of order $n$ and minimum degree $n-1$ that  have no Hamiltonian cycle. As a partial solution of this problem,  \cite{[22]} in his excellent paper proved a structural theorem on the extremal digraphs. An analogous problem  for the Meyniel theorem  was considered by  \cite{[8]}, proving a structural theorem on  the strong non-Hamiltonian digraphs  $D$ of order $n$, with the condition that $d(x)+d(y)\geq 2n-2$ for every pair of non-adjacent  distinct vertices $x,y$. This improves the corresponding   structural theorem of Thomassen.  \cite{[8]} also  proved that: if $m$ is the length of longest cycle in $D$, then $D$ contains cycles of all lengths $k=2,3,\ldots ,m$.  \cite{[22]}  and  \cite{[9]}  described all the extremal digraphs for the Nash-Williams theorem, respectively, when the order of a digraph is odd and   when the order of a digraph is even. 
Here we combine they in the following theorem . \\

  \noindent\textbf{Theorem 1.5} (\cite{[22]} and \cite{[9]}). {\it Let $D$ be a digraph of order $n\geq 4$ with minimum degree $n-1$. If for every vertex $x$ of $D$, $d^+(x)\geq n/2-1$ and $d^-(x)\geq n/2-1$, then $D$ is Hamiltonian, unless some exceptions, which completely are characterized.}\\
 
  \cite{[15]} relaxed the condition of the Ghouila-Houri theorem by proving the following theorem.\\
 
 \noindent\textbf{Theorem 1.6} (\cite{[15]}). {\it Let $D$ be a strong digraph of order $n\geq 2$. If $n-1$ vertices of $D$ have degrees at least $n$ and the remaining vertex has degree at least $n-1$,  then $D$ is Hamiltonian.}\\ 
 
  Note that Theorem 1.6 is an immediate consequence of Theorem 1.4.  
 \cite{[15]} for any $n\geq 5$ presented two examples of non-Hamiltonian strong digraphs of order $n$ such that: 
  (i) In the first example, $n-2$ vertices have degrees equal to  $n+1$ and the other two vertices have degrees equal to $n-1$. 
   (ii) In the second example, $n-1$ vertices have degrees at least $n$ and the remaining vertex has degree equal to $n-2$.\\ 

\textbf{Remark 1.} It is worth to mention that \cite{[22]} constructed a strong  non-Hamiltonian digraph of order $n$ with only two vertices of degree $n-1$ and all other $n-2$ vertices have degrees at least  $(3n-5)/2$.\\
   
 \cite{[28]}  reduced 
the lower bound in Theorem 1.3 by 1, and proved that the conclusion still holds, with only a few exceptional cases that can be clearly characterized.
 \cite{[11]} announced  that the following theorem is holds.\\

 \noindent\textbf{Theorem 1.7} (\cite{[11]}). {\it Let $D$ be a 2-strong  digraph of order $n\geq 9$  such that its $n-1$ vertices have degrees at least  $n$ and the remaining  vertex has degree at least $n-4$.
 Then $D$ is Hamiltonian.}\\

The proof of Theorem 1.7 has never been published. G. Gutin suggested  me to publish the proof of this theorem anywhere. Recently,  \cite{[12]} presented a new proof of the first part of Theorem 1.7, by proving the following:\\

\noindent\textbf{Theorem 1.8} (\cite{[12]}).  {\it Let $D$ be a 2-strong  digraph of order $n\geq 9$  such that its $n-1$ vertices have degrees at least  $n$ and the remaining  vertex $z$ has degree at least $n-4$. If $D$ contains a cycle of length $n-2$ through $z$,
 then $D$ is Hamiltonian.} \\

\cite{[12]} also proposed the following conjecture.

\textbf{Conjecture 1}. {\it Let $D$ be a 2-strong  digraph of order $n$.     Suppose that $n-1$ vertices of $D$ have degrees at least $n+k$ and the remaining vertex has degree at least $n-k-4$, where $k\geq 0$ is an integer. Then $D$ is Hamiltonian.} \\

Note that,  for $k=0$ this conjecture is Theorem 1.7. By inspecting the proof of Theorem 1.8 and the handwritten proof of Theorem 1.7, by the similar arguments we settled  Conjecture 1 by proving the following theorem.\\

\noindent\textbf{Theorem 1.9}.  {\it Let $D$ be a 2-strong  digraph of order $n\geq 9$. If $n-1$ vertices of $D$ have degrees at least $n+k$ and the remaining vertex $z$ has degree at least $n-k-4$,   where $k\geq 0$ is an integer,
then $D$ is Hamiltonian.}\\ 

 \cite{[13]} presented  the proof of the first part of Conjecture 1 for any $k\geq 1$, which we formulated as Theorem 3.6 (Section 3). 
The goal of this article to present the   complete proof of the second part of the proof of Theorem 1.9  and show that this theorem  
 is best possible in the sense that for $k=0$ there is a 2-strong digraph of order $n=8$ (respectively, $n=9$) with the minimum degree $n-4=4$ (respectively, with the minimum degree $n-5=4$) whose $n-1$ vertices have degrees at least $n$, but it is not Hamiltonian.     To see that the theorem is best possible, it suffices consider the
  digraphs defined in the Examples 1 and 2, see Figure 1. In figures an undirected edge represents two directed arcs of opposite directions.\\ 
 
 \textbf{Example 1}. Let $D_8$ be a digraph of order 8 with vertex set $V(D_8)=\{x_1,x_2,x_3,x_4,y_1,y_2,y_3,z\}$ and arc set $A(D_8)$, which satisfies the following conditions: $D_8\langle \{y_1,y_2,y_3\}\rangle$ is a complete digraph, $x_4\rightarrow \{y_1,y_2,y_3\}\rightarrow x_1$, $x_2\rightarrow \{y_1,y_2,y_3\}\rightarrow x_2$, $D_8$ contains the following 2-cycles and arcs
  $x_i\leftrightarrow x_{i+1}$ for all $i\in [1,3]$,  $x_1\leftrightarrow x_{3}$, $x_3\leftrightarrow z$, $x_4\leftrightarrow x_{2}$, $x_4\rightarrow x_1$, $x_4\rightarrow z$ and  $z\rightarrow x_1$. $A(D_8)$ contains no other arcs.\\

\textbf{Example 2}.  Let $D_9$ be a digraph of order 9 with vertex set $V(D_9)=\{x_1,x_2,x_3,x_4,x_5,y_1,y_2,y_3,z\}$ and arc set $A(D_9)$, which satisfies the following conditions: $D_9\langle \{y_1,y_2,y_3\}\rangle$ is a complete digraph, $x_5\rightarrow \{y_1,y_2,y_3\}\rightarrow x_1$,
  $x_3\rightarrow \{y_1,y_2,y_3\}\rightarrow \{x_1,x_2,x_3$\}, $D_9$ contains the following 2-cycles and arcs
  $x_i\leftrightarrow x_{i+1}$ for all $i\in [1,4]$,  $x_1\leftrightarrow x_{4}$, $x_3\leftrightarrow x_5$,
   $x_4\leftrightarrow x_{2}$, $x_4\leftrightarrow z$, $x_5\rightarrow z$, $z\rightarrow x_1$ and $x_5\rightarrow x_1$. $A(D_9)$ contains no other arcs.\\

Observe that every vertex other than $z$ in $D_8$ (in $D_9$) has degree at least $|V(D_8)|=8$ (at least $|V(D_9)|= 9$) and $d(z)=4$ in both digraphs $D_8$ and $D_9$.
 It is not hard to check that for every  $u\in V(D_8)$  ($u\in V(D_9)$), $D_8-u$ ($D_9-u$) is strong, i.e., $D_8$ and $D_9$ both are 2-strong.
To see this, it suffices to consider a longest cycle in $D_8-u$  (in $D_9-u$) and apply the following well-known proposition.\\

\textbf{Proposition 1} (see Exercise 7.26, \cite{[5]}). Let $D$ be a $k$-strong digraph with $k\geq 1$, let $x$ be a new vertex and $D'$ be a digraph obtained from $D$  and $x$  by adding $k$ arcs from $x$ to distinct vertices of $D$ and $k$ arcs from distinct vertices of $D$ to $x$. Then $D'$ also is $k$-strong.\\

Let $D'_9$ be the digraph obtained from $D_9$ by adding the arcs $x_3x_1$ and $x_5x_2$.

Now we will show that $D'_9$ is not Hamiltonian. Assume that this is not the case.  Let $R$ be an arbitrary Hamiltonian cycle in $D'_9$. Then $R$ necessarily   contains either the arc $x_4z$ or the arc $x_5z$. If
$x_4z\in A(R)$, then it is not difficult to see that either $R[x_4,y_i]=x_4zx_1x_2x_3y_i$  or $R[x_4,y_i]=x_4zx_1x_2x_3x_5y_i$, which is impossible since $N^+(y_i, \{x_1,x_2, \ldots , x_5\})=\{x_1,x_2,x_3\}$. We may therefore assume that $x_5z\in A(R)$. Then necessarily $R$ contains the arc $x_3y_i$ and either the path $x_5zx_1$ or the path $x_5zx_4$. It is easy to check that either $x_2x_3\in A(R)$  or $x_4x_3\in A(R)$. If $x_5zx_4$ is in $R$, then $R[x_5,y_i]=x_5zx_4x_j\ldots x_3y_i$, where $j\in [1,3]$, and if    $x_5zx_1$ is in $R$, then $R[x_5,y_i]$  is one of the following pats: $x_5zx_1x_2x_3y_i$, $x_5zx_1x_2x_4x_3y_i$,  $x_5zx_1x_4x_3y_i$ and $x_5zx_1x_4x_2x_3y_i$,    which is impossible since $N^+(y_i, \{x_1,x_2, \ldots , x_5\})=\{x_1,x_2,x_3\}$. So, in all cases we have a contradiction. Therefore, $D'_9$  is not Hamiltonian, which in turn implies that the digraphs $D_9$,
$D_9+\{(x_3x_1)\}$ and $D_9+\{(x_5x_2)\}$  also are not Hamiltonian. By a similar argument we can show that $D_8$ also is  not Hamiltonian.\\
\vspace{-0.5cm}
\begin{figure*}[!h]
\begin{center}
\includegraphics[width=17cm]{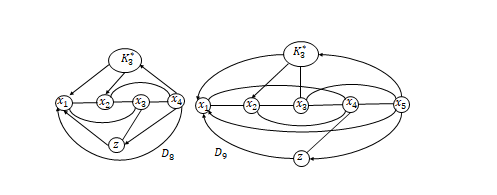}
\vspace{-1cm}
 \caption{The non-Hamiltonian 2-strong digraphs $D_8$ and $D_9$ of order 8 and 9.}
\label{Fig.1}
\end{center}
\end{figure*}

A digraph $D$ is {\it Hamiltonian-connected} if for any pair of distinct vertices $x,y$,  $D$ has a Hamiltonian path from $x$ to $y$. \cite{[21]} proved the following sufficient condition for a digraph to be Hamiltonian-connected. \\

\noindent\textbf{Theorem 1.10} (\cite{[21]}). {\it Let $D$ be a 2-strong  digraph of order $n\geq 3$ such that, for each two non-adjacent distinct vertices $x,y$ we have $d(x)+d(y)\geq 2n+1$. Then for each two distinct vertices $u,v$ with $d^+(u)+d^-(v)\geq n+1$ there is a Hamiltonian $(u,v)$-path.}\\ 
 
  Let $D$ be a   digraph of order $n\geq 3$ and let $u$ and $v$  be two distinct vertices in $V(D)$. Following  \cite{[21]}, we define a new digraph $H_D(u,v)$ as follows:
  $$V(H_D(u,v))=V(D-\{u,v\})\cup \{z\} \,\, (z \,\, \hbox{a new vertex}),$$  
$$A(H_D(u,v))= A(D-\{u,v\})\cup 
\{zy\, |\, y\in N^+_{D-v}(u)\} \cup \{yz\, |\, y\in N^-_{D-u}(v)\}.$$

Now, using Theorem 1.9, we will prove the following theorem, which is an analogue of the Overbeck-Larisch theorem.\\

\noindent\textbf{Theorem 1.11.} {\it Let $D$ be a 3-strong  digraph of order $n+1\geq 10$ and let $u$, $v$ be arbitrary two distinct vertices in $D$.  Suppose that $d^+_D(u)+d^-_D(v)\geq n-k-2$ or  $d^+_D(u)+d^-_D(v)\geq n-k-4$ with $uv\notin A(D)$ and for every vertex $x\in V(D)\setminus \{ u,v\}$,  $d_D(x)\geq n+k+2$. Then $D$ has a Hamiltonian $(u,v)$-path. }
   
  \begin{proof} Let $D$ be a 3-strong  digraph of order $n+1\geq 10$ and let $u, v$ be two distinct vertices in $V(D)$. Suppose that $D$ and $u,v$ satisfy the degree conditions of the  theorem. Now we consider the digraph $H:=H_D(u,v)$ of order $n\geq 9$. By an easy computation, we obtain that the minimum degree of $H$ is at least $n-k-4$, and $H$ has $n-1$ vertices  of degrees at least $n+k$. Moreover, we know that $H$ is 2-strong (see \cite{[10]}). Thus, the digraph $H$ satisfies the conditions of Theorem 1.9. Therefore, $H$ is Hamiltonian, which in turn  implies that in $D$ there is a Hamiltonian  
$(u,v)$-path. \end{proof}

There are a number of sufficient  conditions depending on degree or degree sum for Hamiltonicity 
of  bipartite digraphs. Here we combine several of them in the following theorem.\\

\noindent\textbf{Theorem 1.12}. {\it Let  $D$  be a balanced bipartite strong digraph of order $2a\geq 6$. Then $D$ is Hamiltonian provided one of the following holds:

(a) (\cite{[3]}). $d^+(x)+d^-(y)\geq a+2$ for every pair of vertices   $x$,  $y$ such that $x$,  $y$ belong to different partite sets and $xy\notin A(D)$.

(b) (\cite{[4]}).  $d(x)+d(y)\geq 3a$ for every pair of non-adjacent distinct vertices $x,y$.

(c) (\cite{[1]}).  $d(x)+d(y)\geq 3a$ for every pair of   vertices $x, y$ with a common in-neighbour or a common out-neighbour. 

(d) (\cite{[2]}).  $d(x)+d(y)\geq 3a+1$ for every pair of vertices $x,y$ with a common out-neighbour.}\\

All the lower bounds in Theorem 1.12 are the best possible. However,  \cite{[23]} (respectively,\cite{[26]};  \cite{[25]}) reduced  the lower  bound in Theorem 1.12(a) (respectively, 
Theorem 1.12(b); Theorem 1.12(c)) by one, and completely described all non-Hamiltonian   bipartite digraphs, that is the extremal bipartite digraphs for Theorem 1.12(a) (respectively, Theorem 1.12(b); Theorem 1.12(c)). \cite{[24]} reduced the bound by one in Theorem 1.12(d), but it is Hamiltonian whenever $d(x)+d(y)\geq 3a$ for every  pair of distinct vertices $x, y$ 
 with a common out neighbour. Motivated  by Theorems 1.9, 1.12 and Remark 1, it is natural to suggest the following problems.\\

\textbf{Problem 1}. Suppose that $D$ is a $k$-strong balanced bipartite digraph of order $2a\geq 6$. Let $\{x_0,y_0\}$ be a pair of distinct vertices in $V(D)$ such that $d(x_0)+d(y_0)\geq 3a-l$, where $l\geq 1$ is an integer. Find the minimum value of $k$ and the maximum value of $l$ such that $D$ is Hamiltonian provided one of the following holds:

(i) $x_0$ and $y_0$ are not adjacent and $d(x)+d(y)\geq 3a$ for every pair $\{x,y\}$  of non-adjacent vertices $x$, $y$ other than $\{x_0,y_0\}$.

(ii) $\{x_0,y_0\}$ is a  pair with a common out-neighbour  and $d(x)+d(y)\geq\  3a$ for every  pair $\{x,y\}$ of vertices $x$, $y$ with a common out-neighbour such that 
$\{x,y\}\not= \{x_0,y_0\}$.\\ 

\textbf{Problem 2}. Suppose that $D$ is a $k$-strong balanced bipartite digraph of order $2a\geq 6$. 
Let $u_0$ and $v_0$ be  two vertices from different partite sets such that $u_0 \nrightarrow v_0$ and 
$d^+(u_0)+d^-(y_0)\geq a+2-l$, where $l\geq 2$ is an integer. Find the minimum value of $k$ and the maximum value of $l$ such that $D$ is Hamiltonian provided that the following holds: $d^+(u)+d^-(v)\geq a+2$ for all vertices $u$ and $v$ from different partite sets such that $\{u,v\}\not=\{u_0,v_0\}$ 
and $u\nrightarrow v$.\\\\

\section{Terminology and notation}

In this paper, we consider finite digraphs without loops and multiple arcs.  For the terminology  not defined in this paper, the  reader is referred to  the book  \cite{[5]}. 
The vertex set and the arc set of a digraph $D$ are    denoted 
  by $V(D)$  and   $A(D)$, respectively.  The {\it order} of $D$ is the number
  of its vertices. 
 For any $x,y\in V(D)$, if $xy\in A(D)$,  we also write $x\rightarrow y$, and say that $x$ {\it dominates} $y$ or $y$ is {\it dominated} by $x$. 
 The notion $x y\notin A(D)$ 
 means that $xy\notin A(D)$. If $x\rightarrow y$ and $y\rightarrow x$ we shall use the notation $x\leftrightarrow y$  ($x\leftrightarrow y$ is called  2-{\it cycle}). If $x\rightarrow y$ and $y\rightarrow z$, we write 
$x\rightarrow y\rightarrow z$.
 Let $A$ and $B$ be two disjoint subsets of $V(D)$.
    The notation  $A\rightarrow B$ means that every
   vertex of $A$  dominates every vertex of $B$. 
 We define 
$A_D(A\rightarrow B)=\{xy\in A(D)\, |\, x\in A, y\in B\}$ and $A_D(A,B)=A_D(A\rightarrow B)\cup A_D(B\rightarrow A)$. 
If $x\in V(D)$
   and $A=\{x\}$ we sometimes write $x$ instead of $\{x\}$. The {\it converse digraph} of $D$ is the digraph obtained from $D$ by reversing the direction of all arcs, and   is denoted by $D^{rev}$.
Let $N_D^+(x)$, $N_D^-(x)$ denote the set of  out-neighbors, respectively the set  of in-neighbors of a vertex $x$ in a digraph $D$.  If $A\subseteq V(D)$, then $N_D^+(x,A)= A \cap N_D^+(x)$ and $N_D^-(x,A)=A\cap N_D^-(x)$. 
The {\it out-degree} of $x$ is $d_D^+(x)=|N_D^+(x)|$ and $d_D^-(x)=|N_D^-(x)|$ is the {\it in-degree} of $x$. Similarly, $d_D^+(x,A)=|N_D^+(x,A)|$ and $d_D^-(x,A)=|N_D^-(x,A)|$. The {\it degree} of the vertex $x$ in $D$ is defined as $d_D(x)=d_D^+(x)+d_D^-(x)$ (similarly, $d_D(x,A)=d_D^+(x,A)+d_D^-(x,A)$). We omit the subscript if the digraph is clear from the context. The subdigraph of $D$ induced by a subset $A$ of $V(D)$ is denoted by $D\langle A\rangle$ and $D-A$ is the subdigraph induced by $V(D)\setminus A$, i.e. $D-A=D\langle V(D)\setminus A\rangle$.
For integers $a$ and $b$, $a\leq b$, let $[a,b]$  denote the set $\{x_a,x_{a+1},\ldots , x_b\}$. If $j<i$, then $\{x_i,\ldots , x_j\}=\emptyset$. 
A path is a digraph with vertex set $\{x_1,x_2,\ldots , x_k\}$ and arc set $\{x_1x_2,x_2x_3,\ldots , x_{k-1}x_k\}$, and is denoted by $x_1x_2\ldots x_k$. This is also called an $(x_1,x_k)$--path or a path from $x_1$ to $x_k$. If we add the arc $x_kx_1$  to the above, we obtain a cycle $x_1x_2\ldots x_kx_1$.
 The {\it length} of a cycle or a path is the number of its arcs.
 If a digraph $D$  contains a path from a vertex $x$ to a vertex $y$ we say that $y$ is {\it reachable} from  $x$ in $D$. In particular, $x$ is  reachable from  itself. 
If $P$ is a path containing a subpath from $x$ to $y$, we let $P[x,y]$ denote that subpath. Similarly,  if  $C$ is a cycle containing vertices $x$ and $y$, $C[x,y]$ denotes the subpath of $C$ from $x$ to $y$.   For a cycle $C$,  a $C$-bypass is an $(x,y)$-path $P$ of length at least two  such that   $V(P)\cap V(C)=\{x,y\}$. The {\em flight} of $C$-bypass $P$ respect to $C$ is  $|V(C[x,y])|-2$.

For integers $a$ and $b$, $a\leq b$, let $[a,b]$  denote the set of
all integers, which are not less than $a$ and are not greater than
$b$.

The path (respectively, the cycle) consisting of the distinct vertices $x_1,x_2,\ldots ,x_m$ ($m\geq 2 $) and the arcs $x_ix_{i+1}$, $i\in [1,m-1]$  (respectively, $x_ix_{i+1}$, $i\in [1,m-1]$, and $x_mx_1$), is denoted  by $x_1x_2\cdots x_m$ (respectively, $x_1x_2\cdots x_mx_1$).
We say that $x_1x_2\cdots x_m$ is a path from $x_1$ to $x_m$ or is an $(x_1,x_m)$-{\it path}. Let $x$ and $y$ be two distinct vertices of a digraph $D$. Cycle that passing through $x$ and $y$ in $D$, we denote by $C(x,y)$.
 By $C_m(x)$ (respectively, $C(x)$) we denote   a cycle in $D$ of length $m$ through $x$ (respectively, a cycle through $x$). Similarly, we denote by $C_k$ a cycle of length $k$. By $K_n^*$ is denoted the complete digraph of order $n$. Let $D$ be a digraph of order $n$. If $E$ is a set of arcs in $K_n^*$, then we denote by $D+E$ the digraph  obtained from $D$ by adding all arcs of $E$. 
A digraph $D$ is {\it strongly connected} (or, just, {\it strong}) if there exists a path from $x$ to $y$ and a path from $y$ to $x$ for every pair of distinct vertices $x,y$.     A digraph $D$ is {\it $k$-strongly} connected (or {\it $k$-strong}), where $k\geq 1$, if $|V(D)|\geq k+1$ and $D- A$ is strongly connected for any subset $A\subset V(D)$ of at most $k-1$ vertices. 
Two distinct vertices $x$ and $y$ are {\it adjacent} if $xy\in A(D)$ or $yx\in A(D) $ (or both). We will use {\it 
 the principle of digraph duality}: Let $D$ be a digraph, then $D$ contains a subdigraph $H$ if and only if $D^{rev}$ contains the  subdigraph $H^{rev}$.

\section{Preliminaries}

 In our proofs we extensively will use  the following well-known simple lemmas. 
    
 \noindent\textbf{Lemma 3.1} (\cite{[17]}). {\it Let $D$ be a digraph of order $n\geq 3$ containing a cycle $C_m$, $m\in [2,n-1]$. Let $x$ be a vertex not contained in this cycle. If $d(x,V(C))\geq m+1$, then $D$ contains a cycle $C_k$ for every $k\in [2,m+1]$.}\\
  
  The next lemma is a slight  modification of a lemma by  \cite{[7]} it is very useful and will be used extensively throughout this paper.\\ 
  
 \noindent\textbf{Lemma 3.2.}   Let $D$ be a digraph of order $n\geq 3$ containing a path $P:=x_1x_2\ldots x_m$, $m\in [1,n-1]$.
 Let $x$ be a vertex not contained in this path. If one of the following condition holds: (i) $d(x,V(P))\geq m+2$, (ii) $d(x,V(P))\geq m+1$ and
 $x\nrightarrow x_1$ or $x_m\nrightarrow x$, 
 (iii) $d(x,V(P))\geq m$,  $x\nrightarrow x_1$ and $x_m\nrightarrow x$, 
 then there is an $i\in [1,m-1]$ such that $x_i\rightarrow x\rightarrow x_{i+1}$, i.e.,  $D$ contains a path $x_1x_2\ldots x_ixx_{i+1}\ldots x_m$ of length $m$ (we say that $x$ can be inserted into $P$).\\ 
 
 We note that in the above Lemma 3.2 as well as throughout the whole paper we allow paths  of length 0, i.e., paths that have exactly one vertex. 
Using Lemma 3.2, it is not difficult to prove  the following lemma.\\ 
 
\noindent\textbf{Lemma 3.3.} {\it Let $D$ be a digraph of order $n\geq 4$. Suppose that    $P:=x_1x_2\ldots x_m$, $m\in [2,n-2]$, is a longest path from $x_1$ to $x_m$ in $D$ and $V(D)\setminus V(P)$ contains two distinct vertices $y_1$, $y_2$ such that $d(y_1,V(P))=d(y_2,V(P))=m+1$. 
If in subdigraph $D\langle V(D)\setminus V(P)\rangle$  there exists a path from $y_1$ to $y_2$ and a path from $y_2$ to $y_1$,  then there is an integer $l\in [1, m]$ such that for every $i\in [1,2]$
 $$
O(y_i,V(P))=\{x_1,x_2,\ldots ,x_l\} \quad \hbox{and} \quad I(y_i,V(P))=\{x_l,x_{l+1},\ldots , x_m\}.
$$ }

 \noindent\textbf{Theorem 3.4} (\cite{[10]}). {\it  Let $D$ be a strong digraph of order $n\geq 2$. Suppose that $d(x)+d(y)\geq 2n-1$ for all pairs of  non-adjacent  vertices  $x, y\in V(D)\setminus \{z\}$, where $z$ is an arbitrary fixed vertex in $V(D)$. Then $D$  contains a cycle of length is at least $n-1$.}\\
  
  From Theorem 3.4 it follows that  the following corollary is true.\\   

  \noindent\textbf{Corollary 1}. (\cite{[10]}). Let $D$ be a strong digraph of order $n\geq 2$. Suppose that $n-1$ vertices of $D$ have degrees at least $n$.  Then $D$ either is Hamiltonian or contains a cycle of length $n-1$ (in fact $D$ has a cycle that contains all the vertices with degree at least $n$).\\

\noindent\textbf{Lemma 3.5} (\cite{[12]}). {\it Let $D$ be a digraph of order $n\geq 4$ such that for any vertex
$x\in V(D)\setminus \{z\}$, $d(x)\geq n$, where $z$ is an arbitrary fixed  vertex in $V(D)$. Moreover,  $d(z)\leq n-2$. Suppose that $C_m(z)=x_1x_2\ldots x_mx_1$,  $m\leq n-1$, is a  cycle of length $m$ through $z$ and $C_m(z)$ has an $(x_i,x_j)$-bypass such that $z\notin V(C_m(z)[x_{i+1},x_{j-1}])$. Then $D$ has a cycle, say $Q$, of length at least $m+1$
such that $V(C_m(z))\subset V(Q)$.}\\

\noindent\textbf{Theorem 3.6} (\cite{[13]}).  {\it Let $D$ be a 2-strong  digraph of order $n\geq 9$  such that $n-1$ vertices of $D$ have degrees at least $ n+k$  and the remaining vertex $z$ has degree at least $n-k-4$, where $k\geq 0$ is an integer. If  the length of a longest  cycle through $z$ is at least $n-k-2$, then $D$ is Hamiltonian.}

 \section{Proof of Theorem 1.9}
 
\noindent\textbf{Theorem 1.9.}  {\it Let $D$ be a 2-strong  digraph of order $n\geq 9$. If $n-1$ vertices of $D$ have degrees at least $n+k$ and the remaining vertex $z$ has degree at least $n-k-4$,   where $k\geq 0$ is an integer,
then $D$ is Hamiltonian.} 

 \begin{proof} By contradiction, suppose that $D$ is not Hamiltonian. Then from Theorem 3.6 it follows that $D$ has no $C(z)$-cycle of length greater than $n-k-3$.  By Corollary 1, $D$ contains a cycle of length $n-1$. Let $C_{n-1}:=x_1x_2\ldots x_{n-1}x_1$ be an arbitrary cycle in $D$. By Lemma 3.1, $z\notin V(C_{n-1})$. 
Since $D$ is 2-strong, there are two distinct vertices, say  $x_1$ and  $x_{n-d-1}$, such that $x_{n-d-1}\rightarrow z\rightarrow x_1 $ and $d(z,\{x_{n-d},x_{n-d+1}, \ldots , x_{n-1}\})=0$.
 Without loss of generality, assume that the flight $d:=|\{x_{n-d}, x_{n-d+1},\ldots , x_{n-1}\}|$
    of $z$ respect to $C_{n-1}$ is smallest possible over  all the cycles of length $n-1$ in $D$.

 \vspace{0.5cm}
\begin{figure*}[!h]
\begin{center}
\includegraphics[width=14cm]{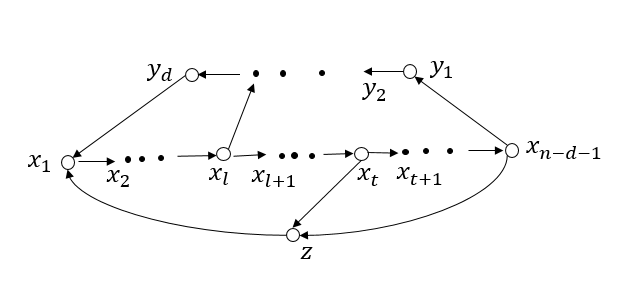}
\vspace{-1cm}
 \caption{The cycles $C_{n-1}=x_1x_2\ldots x_{n-d-1}y_1y_2\ldots y_dx_1$  and $C_{n-d}(z)=x_1x_2\ldots x_{n-d-1}zx_1$ in $D$.}
\label{Fig.2}
\end{center}
\end{figure*}

 For any $i\in [1,d]$, let  $y_i=x_{n-d-1+i}$ and $Y=\{y_1, y_2, \ldots , y_d\}$. Note that $y_1y_2\ldots y_d$ is a path in  $D\langle Y\rangle$. Since $z$ cannot be inserted into $C_{n-1}$, using Lemma 3.2, we obtain $n-k-4\leq d(z)\leq n-d$. 
 Hence, $d\leq k+4$. On the other hand, $n-d\leq n-k-3$, i.e., $d\geq k+3$, since $zx_1x_2\ldots x_{n-d-1}z$ is a $C(z)$-cycle of length $n-d$.
  From now on, by $P$ we denote the path $x_1x_2\ldots x_{n-d-1}$ (see Figure 2).
 In order to prove the theorem, it is convenient for the digraph $D$ and the path $P$ to prove the following  Claims 1-4.\\

\textbf{Claim 1.} {\em Suppose that $D\langle  Y\rangle$ is strong and   each vertex $y_j$ of $Y$ cannot be inserted into $P$.

If $d^+(x_i,Y)\geq 1$ with $i\in [1,n-d-2]$, then $A(Y\rightarrow \{x_{i+1},x_{i+2}, \ldots , x_{n-d-1}\})=\emptyset$}. 

\begin{proof} By contradiction,  suppose that there are vertices $x_s, x_q$ with $1\leq s<q\leq n-d-1$ and $u,v\in Y$ such that $x_s\rightarrow u$, $v\rightarrow x_q$.  Since $D\langle  Y\rangle$ is strong, it contains a $(u,v)$-path, and let $Q$ be  such a longest path.  We may assume that 
$A(Y,\{x_{s+1},\ldots , x_{q-1}\})=\emptyset$. Since $D\langle  Y\rangle$ is strong and every vertex $y_j$ cannot be inserted  into $P$,  using the fact that $D$ has no $C(z)$-cycle  of length at least $n-k-2$, 
   we obtain that $q-s\geq 2$. We now extend the path
 $x_qx_{q+1}\ldots x_{n-d-1}zx_1x_{2}\ldots x_{s}$ with vertices $x_{s+1},x_{s+2},\ldots , x_{q-1}$  as much as possible. Then some vertices $z_1,z_2, \ldots , z_m\in \{x_{s+1},x_{s+2},\ldots , x_{q-1}\}$, where $0\leq m\leq q-s-1$, are not on the obtained extended path, say $R$.  We consider the cases $m\geq 1$ and $m=0$ separately.
 
 Assume first that $m\geq 1$. Since every vertex  $y_j$ cannot be inserted into  $P$  and $d(y_j,\{z,x_{s+1},x_{s+2},\\\ldots ,x_{q-1}\})=0$, using Lemma 3.2(i),  we obtain 
 $$
  n+k\leq d(y_j)=d(y_j,Y)+d(y_j,\{x_1,x_2,\ldots , x_s\})+d(y_j,\{x_q,x_{q+1},\ldots , x_{n-d-1}\}) $$ $$\leq  2d-2+(s+1)+(n-d-1-q+2)=n+s+d-q  \quad \hbox{and}
  $$ 
$$
 n+k\leq d(z_i)=d(z_i,V(R))+d(z_i, \{z_1,z_2,\ldots , z_m\})\leq |V(R)|+1+2m-2 $$ $$= n-d-m+1+2m-2=n+m-d-1.
$$ 
Therefore,
 $$ 
 2n+2k\leq d(z_i)+d(y_j)\leq n+m-d-1+n+s+d-q= 2n+m-1+s-q$$ 
$$
\leq 2n-1+q-s-1+s-q=2n-2,
 $$ 
 which is a contradiction since $k\geq 0$.
 
   Assume next that $m=0$. This means that $D$ contains an $(x_q,x_s)$-path   with vertex set  $\{z\}\cup V(P)$. 
   This and the fact that $D$ contains no cycle of length at least $n-k-2$ through $z$ imply that $d=k+4$, $|V(Q)|=1$, i.e.,  $u=v$,  and
 $
A(x_s\rightarrow Y\setminus \{u\})=A(Y\setminus \{u\}\rightarrow x_q)=\emptyset 
$. 
Since any vertex of $Y$ cannot be inserted into $P$, using  Lemma 3.2(ii),  for each $y\in Y\setminus \{u\}$ we obtain  
 $$
 n+k\leq d(y)=d(y,Y)+d(y,\{x_1,x_2,\ldots , x_s\})+
 d(y,\{x_q,x_{q+1},\ldots , x_{n-k-5}\})$$ $$\leq 2k+6+s+n-k-5-q+1= n+k+2-(q-s).
 $$ 
 This  means that all the inequalities  used in the last expression are actually equalities, i.e.,   $q-s= 2$, $d(y,Y)=2k+6$, i.e.,
 $D\langle Y\rangle$ is a complete digraph, and 
$$
d(y,\{x_1,x_2,\ldots , x_s\})=s,\,\,  d(y,\{x_q,x_{q+1},\ldots , x_{n-k-5}\})=n-k-q-4.$$ 

Again using Lemma 3.2(ii),  from the last two equalities  and $
A(x_s\rightarrow Y\setminus \{u\})=A(Y\setminus \{u\}\rightarrow x_q)=\emptyset$ we obtain    
$x_{n-k-5}\rightarrow Y\setminus \{u\}\rightarrow x_1$. 
We claim that $x_{s+1}$  can be   inserted into $x_1x_2\ldots  x_s$ or $x_qx_{q+1}\ldots x_{n-k-5}$. Assume that this is not the case. Then by Lemma 3.2(i), 
$$
 n+k\leq d(x_{s+1})=d(x_{s+1},\{x_1,x_2,\ldots , x_s\})+
 d(x_{s+1},\{x_q=x_{s+2},x_{s+3},\ldots , x_{n-k-5}\})$$
$$+d(x_{s+1},\{z\})\leq s+1+n-k-5-s-1+1+2=n-k-(q-s)=n-k-2,
 $$  
which is a contradiction. This contradiction shows that there is either an $(x_1,x_s)$-path, say $R_1$, with vertex set $\{x_1,x_2,\ldots , x_s,x_{s+1}\}$  or an $(x_q,x_{n-k-5})$-path, say $R_2$,    with vertex set $\{x_{s+1},x_{s+2},\ldots , x_{n-k-5}\}$. Let $H$ be a Hamiltonian path in
$D\langle Y\setminus \{u\}\rangle$. We know that  $d(z,V(H))=0$, $|V(H)|=k+3$ and $x_{n-k-5}\rightarrow Y\setminus \{u\}\rightarrow x_1$. Therefore, $F_1:=x_1R_1ux_q\ldots x_{n-k-5}Hx_1$ or  $F_2:=x_1\ldots x_suR_2x_{n-k-5}\\Hx_1$, is a cycle of length $n-1$. We have that
  the flight of $z$ respect to $F_1$ (or $F_2$) is equal to $k+3$,
 which contradicts the minimality of $d=k+4$ and the choice of the  cycle $C_{n-1}$ of length $n-1$. This completes the proof of the claim.  \end{proof}
 
 \textbf{Claim 2.} {\it  If
  $x_j\rightarrow z$ with $j\in [1,n-d-2]$, then $A(z\rightarrow \{x_{j+1},x_{j+2},\ldots , x_{n-d-1}\})=\emptyset$.}
  
 \begin{proof}  By contradiction,  suppose that $x_j\rightarrow z$ with $j\in [1,n-d-2]$ and $z\rightarrow x_l$ with $l\in [j+1,n-d-1]$. We may assume that $d(z,\{x_{j+1},\ldots , x_{l-1}\})=0$.  Since $D$ contains no  $C(z)$-cycle of length at least $n-k-2$ and $C_{n-l+j+1}(z):=x_1x_2\ldots x_jzx_l\ldots x_{n-d-1}y_1y_2\ldots y_dx_1$, it follows that $l\geq j+k+4$.   Then, since $z$ cannot be inserted into $P$, by Lemma 3.2(i), we have
 $$
 n-k-4\leq d(z)=d(z,\{x_1,x_2,\ldots , x_j\})+d(z,\{x_l,x_{l+1},\ldots , x_{n-d-1}\})$$ $$\leq (j+1)+(n-d-1-l+2)=n+2+j-d-l$$ $$\leq n+2+(l-k-4)-d-l=n-k-2-d,
 $$
 i.e., $d\leq 2$, which contradicts that $d\geq k+3$. Claim 2 is proved. \end{proof} 
 
 Since $D$ is 2-strong, we have $d^-(z)\geq 2$ and $d^+(z)\geq 2$. From this and Claim 2 it follows that there exists an integer $t\in [2,n-d-2]$ such that $x_t\rightarrow z$ and
 $$
 d^-(z,\{x_1,x_2,\ldots , x_{t-1}\})= d^+(z,\{x_{t+1},x_{t+2},\ldots , x_{n-d-1}\})=0. \eqno (1)
 $$
 From (1) and $d(z)\geq n-k-4$ it follows that if $d=k+4$, then $n-d-1=n-k-5$ and
 $$
N^+(z)=\{x_1,x_2,\ldots , x_{t}\} \quad \hbox{and} \quad N^-(z)=\{x_{t},x_{t+1},\ldots , x_{n-k-5}\}. \eqno (2)
 $$
 
\textbf{Claim 3.}  {\it Suppose that there is an integer $l\in [2,n-d-2]$ such that
  $$
A(\{x_1,x_2, \ldots , x_{l-1}\} \rightarrow Y)=A(Y\rightarrow\{x_{l+1},x_{l+2}, \ldots , x_{n-d-1}\}) = \emptyset.
$$ 
   Then for every  $j\in [2,n-d-2]$, 
 $$
 A(\{x_1,x_2, \ldots , x_{j-1}\}\rightarrow\{x_{j+1},x_{j+2}, \ldots , x_{n-d-1}\}) \not= \emptyset.    
 $$}
 \begin{proof} Suppose, on the contrary, that for some $j\in [2,n-d-2]$, 
 $
 A(\{x_1,x_2, \ldots , x_{j-1}\}\rightarrow\{x_{j+1},x_{j+2}, \\\ldots , x_{n-d-1}\}) = \emptyset    
 $. Without loss of generality, we may assume that $j\leq l$. If  $d^-(z,\{x_1,x_2,\ldots ,  x_{j-1}\})=0$, then by the suppositions of the claim, we have 
$$
A(\{x_1,x_2, \ldots , x_{j-1}\}\rightarrow Y\cup \{z, x_{j+1},x_{j+2}, \ldots , x_{n-d-1}\}) = \emptyset.
$$ 
 If $d^-(z,\{x_1,x_2,\ldots , x_{j-1}\})\geq 1$, then by Claim 2,  $d^+(z,\{x_{j+1},x_{j+2},\ldots , x_{n-d-1}\})=0$. This together with  the supposition of the claim implies that 
 $$
A(\{z,x_1,x_2, \ldots , x_{j-1}\} \rightarrow Y\cup \{x_{j+1},x_{j+2}, \ldots , x_{n-d-1}\}) = \emptyset.
$$
  Thus, in both cases, $D-x_j$ is not strong, which is a contradiction. Claim 3 is proved.  \end{proof}
  
\textbf{Claim 4}. {\em Any vertex $y_j$ with $j\in [1,d]$ cannot be inserted into $P$.}

\begin{proof} By contradiction, suppose that there is a vertex $y_p$ with  $p\in [1,d]$ and an integer $s\in [1,n-d-2]$ such that 
$x_s\rightarrow y_p\rightarrow x_{s+1}$.  Then $R(z):=x_1x_2\ldots x_sy_px_{s+1}\ldots x_{n-d-1}zx_1$ is a cycle of length $n-d+1$. Since $D$ contains no  $C(z)$-cycle of length at least $n-k-2$, it follows that $n-d+1\leq n-k-3$, i.e., $d\geq k+4$. Therefore, $d=k+4$ since $d\leq k+4$.  It is easy to see that any vertex $y_i$ other than $y_p$ cannot be inserted into $P$. Note that (2)  holds since $d=k+4$. We will consider the cases $p\in [2,k+3]$ and $p=1$ separately. Note that if $p=k+4$, then in the converse digraph of $D$ we have case $p=1$. \\ 

\textbf{Case 1.}  $p\in [2,k+3]$.

If $y_{p-1}\rightarrow y_{p+1}$, then  the cycle $x_1x_2\ldots x_sy_px_{s+1}\ldots x_{n-k-5}y_1\ldots y_{p-1}y_{p+1}\ldots y_{k+4}x_1$ is a cycle of length $n-1$ and  the flight of  $z$ respect to this cycle is equal to $k+3$, which is a contradiction. 

 We may therefore assume that $y_{p-1}\nrightarrow y_{p+1}$.
 Since  both $y_{p-1}$ and $y_{p+1}$ cannot be inserted into $R(z)$, using Lemma 3.2(i), we obtain  $d(y_{p-1},V(R(z)))\leq n-k-3$ and $d(y_{p+1}, V(R(z)))\leq n-k-3$. These together with $d(y_{p-1})\geq n+k$ and $d(y_{p+1})\geq n+k$ imply that $d(y_{p-1}, Y\setminus \{y_p\})\geq 2k+3$ and $d(y_{p+1}, Y\setminus \{y_p\})\geq 2k+3$.
Hence, it is easy to see that $y_{p+1}\rightarrow y_{p-1}$ and $d^+(x_s, Y\setminus \{y_p\}) =d^-(x_{s+1}, Y\setminus \{y_p\})=0$ (for otherwise $D$ contains a $C(z)$-cycle of length at least $n-k-2$,  a contradiction). Since every vertex of $Y\setminus \{y_p\}$ cannot be extended into $P$, using Lemma 3.2 and the last equalities, we obtain that if $u\in \{y_{p-1}, y_{p-1}\}$, then 
$$
n+k\leq d(u)=d(u,Y)+ d(u,\{x_1,x_2,\ldots , x_{s}\})+d(u,\{x_{s+1},x_{s+2},\ldots , x_{n-k-5}\})$$ $$\leq 2k+5+s+(n-5-k-s)=n+k.
$$
From this, in particular, we have $d(u,Y)=2k+5$, 
$d(u,\{x_1,x_2,\ldots , x_s\})= s$ and $d(u,\{x_{s+1},x_{s+2},\\\ldots , x_{n-k-5}\})=n-k-5-s$.
Again using Lemma 3.2(ii), we obtain that $x_{n-k-5}\rightarrow \{y_{p-1}, y_{p+1}\}\rightarrow x_1$. From $d(u, Y)=2k+5$ it follows that $u\leftrightarrow Y\setminus \{y_{p-1}, y_{p+1}\}$ since
$y_{p-1}y_{p+1}\notin A(D)$.
Hence it is not difficult to see that in $D\langle Y\setminus \{y_p\}\rangle$ there is a $(y_{p-1},y_{k+4})$- or $(y_{p-1},y_{p+1})$-Hamiltonian path, say $H$. Thus $x_1x_2\ldots x_sy_px_{s+1}\ldots x_{n-k-5}Hx_1$ is a cycle of length $n-1$ and
 the  flight of  $z$  respect to this cycle is equal to $k+3$,  a contradiction.\\
 
 \textbf{Case 2.}  $p=1$, i.e., $x_s\rightarrow y_1 \rightarrow x_{s+1}$.

 Observe that $d^-(x_{s+1},\{y_2,y_3,\ldots ,y_{k+4}\})=0$ and  $R(z)$ is a longest cycle  through $z$ in $D$, which has length $n-k-3$. For Case 2 we will prove the following proposition.\\
 
 \textbf{Proposition 2.} Suppose that for $j$, $j\in [2,k+4]$, in $Q:=D\langle \{y_2, y_3,\ldots , y_{k+4},x_1\}\rangle$ there is a Hamiltonian $(y_j,x_1)$-path, say $H^j$. Then
 $x_{n-k-5}y_j\notin A(D)$. In particular, $x_{n-k-5}y_2\notin A(D)$.

 \begin{proof} Suppose that the claim is not true,  that is $x_{n-k-5}\rightarrow y_j$ with $j\in [2,k+4]$ and $Q$ has a Hamiltonian $(y_j,x_1)$-path, say $H^j$.   Then  $x_1x_2\ldots x_sy_1x_{s+1}\ldots x_{n-k-5}H^jx_1$ is a cycle of length $n-1$ and the flight of $z$ respect to this cycle is equal to $k+3$, a contradiction. 
 Thus $x_{n-k-5}\nrightarrow y_j$. It is easy to see that $H^2=y_2y_3\ldots y_{k+4}x_1$ is a Hamiltonian path in $Q$. Therefore by the first part of this proposition, 
 $x_{n-k-5}\nrightarrow y _2$. \end{proof} 

To complete the proof of Claim 4, we will consider the cases $x_s\nrightarrow y _2$,
$x_s\rightarrow y_2$ separately.

\textbf{Subcase 2.1.}  $x_sy_2\notin A(D)$.  

We know that    
$y_2x _{s+1}\notin A(D)$ and $x_{n-k-5}y_2\notin A(D)$.
Now, since  $y_2$  cannot be inserted into $P$,  using Lemmas 3.2(ii) and 3.2(iii), we obtain
$$
n+k\leq d(y_2)=d(y_2,Y)+d(y_2,\{x_1,x_2,\ldots , x_s\})$$ $$ +d(y_2,\{x_{s+1},x_{s+2},\ldots , x_{n-k-5}\})\leq 2k+6+s+(n-k-s-6)=n+k.
$$
 This implies that  $d(y_2,Y)=2k+6$, i.e., 
$y_2\leftrightarrow Y\setminus \{y_2\}$, in particular, $y_2\leftrightarrow  y_1$ and $D\langle Y\rangle$ is strong, and 
$$
d(y_2,\{x_1,x_2, \ldots ,\\ x_s\})=s \,\, \hbox{and} \,\, d(y_2,\{x_{s+1},x_{s+2}, \ldots , x_{n-k-5}\})=n-k-s-6. \eqno (3)
$$  
 Thus, for the longest cycle $R(z)$ we have that $V(D)\setminus V(R(z))=\{y_2,y_3,\ldots ,y_{k+4}\}$, $D\langle V(D)\setminus V(R(z))\rangle$  is strong and 
$y_2\leftrightarrow  y_1$.  Therefore by Lemma 3.5, 
$$
 A(\{y_2,y_3,\ldots , y_{k+4}\}\rightarrow \{x_{s+1},x_{s+2}, \ldots , x_{n-k-5}\})=
 A(\{x_1,x_2, \ldots , x_{s}\}\rightarrow \{y_2,y_{3}, \ldots , y_{k+4} \}) = \emptyset.\eqno (4)
$$
This together with $x_{n-k-5}\nrightarrow y_2$ 
and  (3) implies that $y_2$ and $x_{n-k-5}$ are not adjacent and
$$
N^+(y_2, V(P))=\{x_1,x_2, \ldots , x_{s}\} \,\, \hbox{and} \,\,
N^-(y_2,V(P))=\{x_{s+1},x_{s+2}, \ldots , x_{n-k-6}\}.  
$$
By the above arguments, we have that $H^3=y_3y_4\ldots y_{k+4}y_2x_1$ is a $(y_3,x_1)$-Hamiltonian path in $Q$. Therefore by Proposition 1, $x_{n-k-5}\nrightarrow y_3$. This together with (4) implies that $x_{n-k-5}$ and $y_3$ are not adjacent. As for $y_2$, for $y_3$ we obtain that 
$y_3\leftrightarrow Y\setminus \{y_3\}$ and 
$$
N^+(y_3, V(P))=\{x_1,x_2, \ldots , x_{s}\} \,\, \hbox{and} \,\,
N^-(y_3,V(P))=\{x_{s+1},x_{s+2}, \ldots , x_{n-k-6}\}. 
$$
Proceeding in  the same manner, we obtain that $d(x_{n-k-5},\{y_2,y_3, \ldots , y_{k+4}\})=0$, $D\langle Y \rangle$ is a complete digraph and 
$$
\{x_{s+1},x_{s+2}, \ldots , x_{n-k-6}\}\rightarrow Y\setminus \{y_1\} \rightarrow \{x_1,x_2, \ldots , x_{s}\}. \eqno (5)
$$

If $s=n-k-6$, then from (4) and $d(x_{n-k-5},\{y_2,y_3,\ldots , y_{k+4}\})=0$ 
it follows that $A(V(P)\cup \{z\}\rightarrow Y\setminus \{y_1\})=\emptyset$, i.e., 
 $D-y_1$ is not strong, a contradiction.  Therefore, we may assume that $s\leq n-k-7$. 
Let $s=1$. Since $D\langle Y\rangle$ is strong, from (5) it follows that $A(Y\rightarrow \{x_{3},x_4,\ldots , x_{n-k-5}\})=\emptyset$.
(for otherwise, $y_1\rightarrow x_i$ with $i\in [3,n-k-5]$ and $C_{n-k-2}(z) =x_1x_2\ldots x_{i-1}y_2y_1x_i\ldots x_{n-k-5}zx_1$, a contradiction).
If $d^+(x_1,\{x_3,x_4,\ldots , x_{n-k-5}\})=0$, then $A(\{x_1\}\cup Y\rightarrow \{z,x_3,x_4, \ldots , x_{n-k-5}\})=\emptyset$, i.e.,  $D-x_2$ is not strong, a contradiction. So, we can assume that
 for some $b\in [3,n-k-5]$, $x_1\rightarrow x_b$. By (5) and (2), respectively, we have $x_{b-1}\rightarrow y_2$ and $z\rightarrow x_2$. Therefore, $C_{n-1}(z):=x_1x_b\ldots x_{n-k-5}zx_2\ldots x_{b-1}y_2y_3\ldots y_{k+4}x_1$, a contradiction. Let finally $2\leq s\leq n-k-7$. It is easy to see that $A(\{x_1,x_2,\ldots , x_{s-1}\}\rightarrow \{x_{s+1},x_{s+2}, \ldots ,x_{n-k-5} \})\not=\emptyset$
(for otherwise, using  the fact that $A(\{x_1,x_2,\ldots , x_{s-1}\}\rightarrow Y)=\emptyset$ (by (5)), Claim 2 and (2),  it is not difficult to show  that  $D-x_s$ is not strong,  a contradiction). 
 Thus, there are integers $a\in [1,s-1]$ and $b\in [s+1,n-k-5]$ such that  $x_a\rightarrow x_b$. 
Then by (4), $y_2\rightarrow x_{a+1}$, and by (2), either $z\rightarrow x_{a+1}$  or $x_{b-1}\rightarrow z$. 
By (4), we also have that $x_{b-1}\rightarrow y_2$ or  $x_{b-1}\rightarrow y_1$ when $b=s+1$. 
Therefore, $C(z)=x_1x_2\ldots x_ax_b\ldots x_{n-k-5}zx_{a+1}\ldots x_{b-1}(y_1\, or \,y_2)y_2y_3\ldots y_{k+4}x_1$ is a cycle of length at least $n-1$ or
 $C_{n-k-2}(z)=x_1x_2\ldots x_ax_b\ldots x_{n-k-5}y_1y_2x_{a+1}\ldots x_{b-1}zx_1$,  respectively, for $z\rightarrow x_{a+1}$ and for $x_{b-1}\rightarrow z$. Thus, for any possible case we have  a contradiction. This completes the discussion  of Subcase 2.1. 
 
\textbf{Subcase 2.2.}  $x_s\rightarrow y_2$.

Using Lemma 3.5 and the fact that $R(z)$  is a longest  cycle of length $n-k-3$  through $z$, we obtain 
$$
A(Y\setminus \{y_1\}\rightarrow \{x_{s+1},x_{s+2},\ldots , x_{n-k-5}\})=\emptyset. \eqno (6)
$$
Since $x_s\rightarrow y_1\rightarrow  x_{s+1}$, it follows that  in $D\langle Y\rangle$ there is no $(y_2,y_1)$-path, i.e., $d^-(y_1,\{y_2,y_3,\ldots ,\\ y_{k+4}\})=0$
 (for otherwise $D$ has a cycle of length  at least $n-k-2$ through $z$, which is a contradiction).
   This implies that for all $i\in [1,k+4]$, $d(y_i,Y)\leq 2k+5$. Recall that 
   $x_{n-k-5}\nrightarrow y_2$ (Proposition 1).
   Therefore, since $y_2$ cannot be inserted into $P$ 
   and $y_2\nrightarrow x_{s+1}$, using Lemma 3.2, we obtain 
 $$
n+k\leq d(y_2)=d(y_2,Y)+d(y_2,\{x_1,x_2,\ldots , x_s\})+d(y_2,\{x_{s+1},x_{s+2},\ldots , x_{n-k-5}\})$$ $$\leq 2k+5+s+1+(n-k-s-6)=n+k.
$$
Therefore,  $y_2\leftrightarrow Y\setminus \{y_1,y_2\}$, in particular, $D\langle Y\setminus \{y_1\}\rangle$ is strong,
$$
d(y_2,\{x_1,x_2,\ldots , x_s\})=s+1 \,\, \hbox{and} \,\, d(y_2,\{x_{s+1},x_{s+2}, \ldots , x_{n-k-5}\})=n-k-s-6. \eqno (7)
$$ 
From (6) and $x_{n-k-5} y_2\notin A(D)$ 
it follows that $y_2$ and $x_{n-k-5}$ are not adjacent. 
Therefore by (7) and (6), $\{x_{s+1},x_{s+2}, \ldots , x_{n-k-6}\}\rightarrow y_2$, and by Lemma 3.2, $y_2\rightarrow x_1$. Note that  $H^3=y_3y_4\ldots y_{k+4}y_2x_1$ is a Hamiltonian $(y_3,x_1)$-path in $Q$. Therefore by Proposition 1, $x_{n-k-5}y_3\notin A(D)$, which together with (6) implies that $y_3$ and $x_{n-k-5}$ are not adjacent.  
Now by the same arguments, as for $y_2$, we obtain that $y_3\leftrightarrow Y\setminus \{y_1,y_3\}$,
$$
d(y_3,\{x_1,x_2,\ldots , x_s\})=s+1 \,\, \hbox{and} \,\,
d(y_3,\{x_{s+1},x_{s+2}, \ldots , x_{n-k-6}\})=n-k-s-6. \eqno (8)
$$ 
Now by (8) and (6), $\{x_{s+1},x_{s+2}, \ldots , x_{n-k-6}\}\rightarrow y_3$. We know that $P_1:=x_1x_2\ldots x_s$ is a longest $(x_1,x_s)$-path in 
$D\langle V(P_1)\cup Y\setminus \{y_1\}\rangle$. Therefore, since $d(y_2,V(P_1))=d(y_3,V(P_1))= s+1$, by Lemma 3.3, there exists an integer $q\in [1,s]$ such that for every $j\in [2,3]$
$$
N^+(y_j,V(P_1))=\{x_1,x_2,\ldots , x_q\} \,\, \hbox{and}\,\, N^-(y_j,V(P_1))=\{x_q,x_{q+1},\ldots , x_s\}.
$$
Proceeding in the same manner, we conclude that $\{x_{s+1},x_{s+2}, \ldots , x_{n-k-6}\}\rightarrow Y\setminus \{y_1\}$, 
 for all $j\in[2,k+d]$,  the vertices $y_j$ and $x_{n-k-5}$ are not adjacent and
$$
N^+(y_j,V(P_1))=\{x_1,x_2,\ldots , x_q\} \,\, \hbox{and}\,\, N^-(y_j,V(P_1))=\{x_q,x_{q+1},\ldots , x_s\}. \eqno (9)
$$
If $q=1$, then $A(\{y_2,y_3,\ldots , y_{k+4}\}\rightarrow \{z,y_1,x_2,x_3,\ldots , x_{n-k-5}\})=\emptyset$, which implies that $D- x_1$ is not strong, a contradiction. Therefore, we may assume that $q\geq 2$, i.e., $q\in [2,s]$. If $x_i\rightarrow y_1$ with $i\in [1,q-1]$ then by (9), $C_n(z)=x_1x_2\ldots x_iy_1y_2\ldots y_{k+4}x_{i+1}x_{i+2}\ldots
 x_{n-k-5}zx_1$, a contradiction. We may therefore assume that $d^-(y_1,\{x_1,x_2,\ldots, x_{q-1}\})=0$. This together with (9) implies that 
$A(\{x_1,x_2,\ldots, x_{q-1}\}\rightarrow Y)=\emptyset$. Since $D$ is 2-strong, the last equality and (2) imply that there are integers $a\in [1,q-1]$ and $b\in [q+1,n-k-5]$ such that $x_a\rightarrow x_b$, for otherwise it is easy to see that $D-x_q$ is not strong. By (9) and (2),  we have $y_{k+4}\rightarrow x_{a+1}$, $x_{b-1}\rightarrow y_2$ and
$z\rightarrow x_{a+1}$ or $x_{b-1}\rightarrow z$. Therefore, if $z\rightarrow x_{a+1}$, then $C_{n-1}(z)=x_1x_2\ldots x_ax_b\ldots x_{n-k-5}
zx_{a+1}\ldots x_{b-1}y_2\ldots y_{k+4}x_1$,  and if $x_{b-1}\rightarrow z$, then 
$C_{n}(z)=x_1x_2\ldots x_ax_b\ldots x_{n-k-5}y_1\ldots y_{k+4}x_{a+1}
\ldots  x_{b-1}zx_1$. 
So, in any case we have a contradiction.  Claim 4 is proved. \end{proof}

For any $j\in [1,d]$, we have 
 $$n+k\leq d(y_j)=d(y_j,V(P))+d(y_j,Y)\leq d(y_j,V(P))+2d-2.$$ 
 From this, $d(y_j,V(P))\geq n+k-2d+2$. On the other hand, by Lemma 3.2 and Claim 4, $d(y_j, V(P))\leq n-d$. Therefore,
$$
 n+k-2d+2\leq d(y_j, V(P))\leq n-d \quad \hbox{and} \quad d+k\leq d(y_j,Y)\leq 2d-2. \eqno (10)
 $$  
We distinguish two cases according to the subdigraph 
$D\langle Y\rangle$ is strong or not.\\

\textbf{Case A.} $D\langle Y\rangle$ is strong.

In this case, by Claim 4, the suppositions of Claim 1 hold. Therefore,  if for some 
$$ i\in [1,n-d-2] \,\, \hbox{and}\,\,
d^+(x_i,Y)\geq 1,\,\, \hbox{then}\,\, A(Y\rightarrow \{x_{i+1},x_{i+2}, \ldots , x_{n-d-1}\})=\emptyset. \eqno (11)
$$
Since $D$ is 2-strong, (11) implies that $d^+(x_1,Y)=d^-(x_{n-d-1},Y)=0$, there exists $l\in [2,n-d-2]$ such that $d^+(x_l,Y)\geq 1$ and
$$
 A(\{x_1,x_2,\ldots , x_{l-1}\}\rightarrow Y)=A(Y\rightarrow
 \{x_{l+1},x_{l+2},\ldots , x_{n-d-1}\})=\emptyset. \eqno (12)
 $$
From this we see that the supposition of Claim 3 holds. Therefore, for all  
  $j\in [2,n-d-2]$,
 $$
 A(\{x_1,x_2, \ldots , x_{j-1}\}\rightarrow\{x_{j+1},x_{j+2}, \ldots , x_{n-d-1}\}) \not= \emptyset.     \eqno(13)
 $$
For Case A, we will  prove the following two claims.\\

 \textbf{Claim 5.} (i) {\em $A(D)$ contains every arc of the forms  $z\rightarrow x_i$ and $x_j\rightarrow z$, where $i\in [1,t]$  and $j\in [t,n-d-1]$, 
    maybe except one when $d=k+3$. (Recall that the definition of $t$ is given immediately after the proof of Claim 2). 
 
 (ii) For every $i\in [1,d]$, $A(D$) contains every arc of the forms $y_i\rightarrow x_q$ and $x_j\rightarrow y_i$ where $q\in [1,l]$ and $j\in [l,n-d-1]$,  
    maybe except one when $d=k+3$ or except two when $d=k+4$.}  
 
\begin{proof} (i) If $d=k+4$, then Claim 5(i) is an immediate consequence of (2). Assume that $d=k+3$. Then by (1), we have  
$$
n-k-4\leq d(z)=d^+(z,\{x_1,x_2,\ldots , x_{t-1}\})+
d(z,\{x_t\})+ d^-(z,\{x_{t+1}, x_{t+2},\ldots , x_{n-k-4}\})$$ $$
\leq t-1+2+ n-k-4-t=n-k-3.$$
Now, it is easy to see that Claim 5(i) is true.

(ii) By (10) and (12) we have 
$$
n+k-2d+2\leq d(y_i,V(P))=d^+(y_i,\{x_1,x_2,\ldots , x_{l-1}\})+
d(y_i,\{x_l\})$$ $$+ d^-(y_i,\{x_{l+1},x_{l+2},\ldots , x_{n-d-1}\})\leq l-1+2+n-d-1-l=n-d.
$$
Now, considering the cases $d=k+3$ and $d=k+4$ separately, it is not difficult to see that Claim 5(ii) also is true. Claim 5 is proved.
\end{proof}

  \textbf{Claim 6}. {\em  Suppose that for some integers $a$ and $b$ with $1\leq a<b-1\leq n-d-2$ we have $x_a\rightarrow x_b$. If 
$D\langle Y\rangle$ is strong and $z\rightarrow x_{a+1}$, then $d^+(x_{b-1},Y)=0$.} 

\begin{proof}  Suppose, on the contrary,  that is $D\langle Y\rangle$ is strong, $z\rightarrow x_{a+1}$ and  $d^+(x_{b-1}, Y)\geq 1$. Let $x_{b-1}\rightarrow y_i$,  where $i\in [1,d]$. Recall that $k+3\leq d\leq k+4$. If $i\in [1,k+3]$, then the cycle $C(z)=x_1x_2\ldots x_ax_b\ldots x_{n-d-1}zx_{a+1}\ldots x_{b-1}y_i\ldots y_dx_1$ has length at least $n-k-2$, which is a contradiction. Therefore, we may assume that  $d^+(x_{b-1},\{y_1,y_2,\ldots , y_{k+3}\})=0$. Then from $d^+(x_{b-1},Y)\geq 1$ it follows that $d=k+4$ and $x_{b-1}\rightarrow y_{k+4}$.  Hence by (11), $A(Y\rightarrow \{x_{b},x_{b+1},\ldots  , x_{n-k-5}\})=\emptyset$. 
Note that for each $i\in [1,k+3]$, $D\langle Y\rangle$ contains a $(y_{k+4},y_i)$-path since  $D\langle Y\rangle$ is strong. Hence it is not difficult to see that if 
 $d^-(x_1, \{y_1,y_2,\ldots , y_{k+3}\})\geq 1$, then $D$ contains  a  $C(z)$-cycle of length at least $n-k-2$, a contradiction. Therefore, we may assume that 
  $d^-(x_1, \{y_1,y_2,\ldots , y_{k+3}\})= 0$. This together with $d^+(x_1,Y)=0$ implies that $d(x_1, \{y_1,y_2,\ldots ,y_{k+3}\})= 0$.  Now using Lemma 3.2, Claim 4, $A(Y\rightarrow \{x_b,x_{b+1},\ldots , x_{n-k-5}\})=\emptyset$ and 
$d^+(x_{b-1},\{y_1,y_2,\ldots , y_{k+3}\})=0$, for any $i\in [1,k+3]$ we obtain, 
$$
n+k\leq d(y_i)=d(y_i, Y)+d(y_i,\{x_2, x_3,\ldots , x_{b-1}\})+d^-(y_i,\{x_b,x_{b+1},\ldots , x_{n-k-5}\})$$ $$ \leq 2k+6+(b-2)+(n-k-5-b+1)=n+k.
$$
This means that all  inequalities which  were used   in  the last expression in fact are equalities, i.e., for any $i\in [1,k+3]$,
$d(y_i, Y)=2k+6$ (i.e., $D\langle Y\rangle$ is a complete digraph), 
and $d(y_i,\{x_2,x_3, \ldots , x_{b-1}\})=b-2$.
Therefore, since any vertex $y_i$ with $i\in [1,k+3]$ cannot be inserted into $P$ (Claim 4), 
$d(y_i,\{x_2,x_3, \ldots ,\\ x_{b-1}\})=b-2$ and 
$x_{b-1}\nrightarrow y_i$, using Lemma 3.2,
we obtain  that  $y_i\rightarrow x_2$.  Hence, if $a\geq 2$, then 
 $C_{n-1}(z)=x_2\ldots x_ax_b\ldots x_{n-k-5}zx_{a+1}\ldots x_{b-1}Lx_2$, where $L$ is a Hamiltonian  $(y_{k+4},y_{k+3})$-path in $D\langle Y\rangle$, a contradiction. Therefore, we may assume that $a=1$. Then $z\rightarrow x_2$ and $C_{n-1}=x_1x_a\ldots x_{n-k-5}y_1y_2\ldots y_{k+3}x_2\ldots x_{b-1}y_{k+4}x_1$ is a cycle of length $n-1$ in $D$. We have  $x_{n-k-5}\rightarrow z$ and $z\rightarrow x_2$, i.e., the flight of $z$ respect to this cycle $C_{n-1}$ is equal to $k+3$, which contradicts that the minimal flight of $z$ respect to on all cycles of length $n-1$ is equal to $d=k+4$.
Claim 6 is proved. \end{proof}\\

Now using the  digraph duality, we  prove that  it suffices to consider only the case $t\geq l$.

 Indeed, assume that $l\geq t+1$ and consider the converse digraph $D^{rev}$ of $D$. Let $V(D^{rev})=\{u_1,u_2,\ldots , u_{n-d-1},v_1,v_2,\ldots ,  v_d,z\}$,
 where $u_i:=x_{n-d-i}$ and $v_j:=y_{d+1-j}$ for all $i\in [1,n-d-1]$ and $j\in [1,d]$, in particular, $x_l=u_{n-d-l}$ and $x_t=u_{n-d-t}$. Let $p:=n-d-l$ and $q:=n-d-t$. Note that $q\geq p+1$ and $\{v_1,v_2, \ldots , v_d\}=Y$.

Observe that from the definitions of $l, t, p$ and $q$ it follows that $d^-_{D^{rev}}(u_p,Y)\geq 1$  $zu_q\in A(D^{rev})$, $d^-_{D^{rev}}(z,\{u_1,u_2,\ldots , u_{q-1}\})=0$  and  $A_{D^{rev}}(\{u_1,u_2,\ldots , u_{p-1}\}\rightarrow Y)=\emptyset$. Now using Claim 5(i), we obtain that
  $d^-_{D^{rev}}(z,\{u_q,u_{q+1}\})\geq 1$  and  $A_{D^{rev}}(\{u_p,u_{p+1}\}\rightarrow Y)\not=\emptyset$ when $d=k+3$ and   $A_{D^{rev}}(\{u_p,u_{p+1}, u_{p+2}\}\rightarrow Y)\not=\emptyset$ when $d=k+4$. Let $u_{t'}z\in A(D^{rev})$, $d^+_{D^{rev}}(u_{l'},Y)\geq 1$ and $t', l'$ are minimal with these properties. It is clear that $t'\in [q,q+1]$ and $l'\in [p,p+2]$.
    We claim that $t'\geq l'$. Assume that this is not the case, i.e., $t'\leq l'-1$. Then it is not difficult to see that $t'\leq l'-1$ is possible when $l'=p+2$ and $t'=p+1=q$.  By Claim 5(ii), $d=k+4$ and $2\leq p=q-1\leq n-k-7$. Therefore, in $D^{rev}$ the following hold: 
    $$
    d_{D^{rev}}(u_{p+1},Y)=0, \,\, \{u_{q+1},u_{q+2},\ldots , u_{n-k-5}\}\rightarrow Y \rightarrow \{u_{1},u_2,\ldots , u_{p}\},
    $$ 
    $$N^+_{D^{rev}}(z)\\=\{u_{1},u_2,\ldots , u_{q}\}\,\, \hbox{and} \,\, 
     N^-_{D^{rev}}(z)=\{u_{q},u_{q+1},\ldots , u_{n-k-5}\}.$$ 
     Since $D^{rev}$ is 2-strong and  
$A_{D^{rev}}(\{u_{1},u_2,\ldots , u_{p}\}\rightarrow \{z\}\cup Y)=\emptyset$, 
it follows that there are  $r\in [1,p]$ and $s\in [p+2,n-k-5]$ such that $u_ru_s\in A(D^{rev})$. 
Taking into account the above observation, it is not difficult to show that if $r\leq p-1$, then $C_n(z)=u_1u_2\ldots u_ru_s\ldots u_{n-k-5}v_1v_2\ldots v_{k+4}u_{r+1}\ldots \\ u_{s-1}zu_1$ 
is a Hamiltonian cycle in $D^{rev}$, and if $s\geq p+3=q+2$, 
then $C_n(z)=u_1u_2\ldots u_ru_s\ldots \\ u_{n-k-5}zu_{r+1}\ldots  u_{s-1}v_1v_2\ldots v_{k+4}u_1$
is a Hamiltonian cycle in $D^{rev}$, which contradicts that $D$ is not Hamiltonian. 
 We may therefore assume that $r=p$ and $s=p+2$.  This means that 
 $A_{D^{rev}}(\{u_1,u_{2},\ldots ,\\ u_{p-1}\}\rightarrow \{u_{p+2},u_{p+3},\ldots ,  u_{n-k-5}\})=\emptyset$. 
Therefore, since   $D^{rev}$ is 2-strong, for some $i\in [1,p-1]$, $u_iu_{p+1}\in  A(D^{rev})$.
Hence, $u_1u_2\ldots u_iu_{p+1}zx_{i+1}\ldots u_pu_{p+2}\ldots u_{n-k-5}v_1v_2\ldots v_{k+4}u_1$ is a Hamiltonian cycle in $D^{rev}$,  a contradiction. 
Therefore,  the case $t\leq l-1$ is equivalent to  the case $t\geq l$.\\ 

Using Lemma 3.1, it is easy to see that the following proposition holds.

\textbf{Proposition 3}. If $k=0$ and a longest $C(z)$-cycle in $D$ has length $n-3$, then $D\langle V(D)\setminus V(C(z))\rangle$ is strong.\\
 
From now on, we assume that $l\leq t$. Note that from (13) it follows that there are $a\in [1,t-1]$ and $b\in [t+1,n-d-1]$ such that $x_a\rightarrow x_b$.\\

\textbf{Subcase A.1.}  
$z\rightarrow x_{a+1}$.

Recall that $a\in [1,t-1]$ and $b\in [t+1,n-d-1]$.
By Claim 6, we have that $d^+(x_{b-1},Y)=0$.\\ 

\textbf{Subcase A.1.1}. $z\rightarrow x_{a+1}$ and $b\geq t+2$. 

Then $b-2\geq t\geq l$. If $x_{b-2}\rightarrow y_i$ with   $i\in [1,2]$, then the cycle $C(z)=x_1x_2\ldots x_ax_b\ldots x_{n-d-1}\\zx_{a+1}\ldots x_{b-2}y_i \ldots y_dx_1$ 
has length at least $n-2$, a contradiction.  We may therefore assume that $d^+(x_{b-2},\{y_1,y_2\})=0$. This together with Claim 6 implies that $A(\{x_{b-2}, x_{b-1}\}\rightarrow \{y_1,y_2\})=\emptyset$. 
Therefore by  Claim 5(ii) and $l\leq t$,  we have that $d=k+4$, in particular, (2) holds. If $b\geq t+3$, then from $d^-(y_1,\{ x_{b-2}, x_{b-1}\})=0$ and Claim 5(ii) it follows that  $x_{b-3}\rightarrow y_1$  and 
 $C_{n-2}(z)=x_1x_2\ldots x_ax_b\ldots x_{n-k-5}zx_{a+1}\ldots x_{b-3}y_1y_2\ldots y_{k+4}x_1$, a contradiction. Therefore,
we may assume that $b=t+2$.  If  $x_{t}\rightarrow y_3$, then  $C_{n-3}(z)=x_1x_2\ldots x_ax_{t+2}\ldots x_{n-k-5}zx_{a+1}\ldots x_{t}y_3\ldots y_{k+4}x_1$ and the subdigraph 
$D\langle V(D)\setminus V(C_{n-3}(z))\rangle = D\langle \{x_{t+1},y_1,y_2\}\rangle $
is not strong since $d^+(x_{t+1},\{y_1,y_2\})=0$. This implies that $n-3\leq n-k-3$ (i.e., $k=0$) since the length of a longest $C(z)$-cycle is at most $n-k-3$. So, we have  a contradiction to Proposition 3. Therefore, $d^-(y_3,\{x_t,x_{t+1}\})=0$. Hence, if $l=t$, then $y_3\rightarrow x_{a+1}$ (Claim 5(ii)) and $C_{n-k-2}(z)=x_1x_2\ldots x_ax_{t+2}\ldots x_{n-k-5}y_1y_2y_3x_{a+1}\ldots x_tzx_1$, a contradiction. 
We may assume that $l\leq t-1$. If $a\leq t-2$, then from $d^-(y_1,\{x_t,x_{t+1}\})=0$ and Claim 5(ii), we have $x_{t-1}\rightarrow y_1$ and $C_{n-2}(z)=x_1x_2\ldots x_ax_{t+2}\ldots x_{n-k-5}zx_{a+1}\ldots x_{t-1}y_1y_2\ldots y_{k+4}x_1$, a contradiction.
We may therefore assume that $a=t-1$. From $l\leq t-1$ and $l\geq 2$ it follows that $t\geq 3$. Thus we have that $a=t-1\geq 2$ and $b=t+2$, which mean that $A(\{x_1,x_2,\ldots ,x_{a-1}=x_{t-2}\}\rightarrow \{x_{t+2},x_{t+3},\ldots , x_{n-k-5}\})=\emptyset$. This together with   (13) implies that  for some 
$i\in [1,t-2]$ and $j\in [t,t+1]$, $x_i\rightarrow x_j$. Recall that $z\rightarrow x_{i+1}$ and  $x_{t+1}\rightarrow z$ because of (2). 
Therefore, $C(z)=x_1x_2\ldots x_ix_{j} x_{t+1}zx_{i+1}\ldots x_{t-1}x_{t+2}\ldots x_{n-k-5}y_1y_2\ldots y_{k+4}x_1$ is cycle of length at least $n-1$,
 a contradiction. This completes the discussion of Sabcase A.1.1.\\ 

\textbf{Subcase A.1.2}. $z\rightarrow x_{a+1}$ and $b= t+1$.
Since $b-1=t$,  $d^+(x_{b-1}, Y)=0$ (Claim 6) and $d^+(x_l, Y)\geq 1$, we have  $d^+(x_t, Y)=0$, $t-1\geq l\geq 2$. 

Assume first that $t+1\leq n-d-2$. Taking into account  Subcase A.1.1 and  $b= t+1$, we may assume that   $A(\{x_1,x_2,\ldots , x_{t-1}\}\rightarrow 
 \{x_{t+2},x_{t+3},\ldots , x_{n-d-1}\})=\emptyset$.
This together with (13) implies that there is  $j\in[t+2,n-d-1]$ such that $x_t\rightarrow x_j$. If 
 $x_{j-1}\rightarrow z$, then $C_n(z)= x_1x_2\ldots x_ax_{t+1}\ldots x_{j-1}zx_{a+1}\ldots x_t\\x_j\ldots x_{n-d-1}y_1y_2\ldots y_dx_1$, a contradiction.
  Therefore, we may assume that $x_{j-1}\nrightarrow z$.
  This together with (2) implies that $d=k+3$. If $j\geq t+3$, then $x_{j-2}\rightarrow z$ (Claim 5(i)) and  $C_{n-1}(z)= x_1x_2\ldots x_ax_{t+1}\ldots x_{j-2}zx_{a+1}\ldots x_tx_j\ldots x_{n-k-4}y_1y_2\ldots y_{k+3}x_1$, a contradiction. Assume that $j=t+2$. Since $d^+(x_t, Y)=0$, $d^+(x_l, Y)\geq 1$, $d=k+3$ and  $l\leq t-1$, by Claim 5(ii) we have $\{x_l,x_{l+1},\ldots ,x_{t-1}\}\rightarrow y_1$. If $a\leq t-2$, then $x_{t-1}\rightarrow y_1$ and $C_{n-1}(z)=x_1x_2\ldots x_ax_{t+1}\ldots x_{n-k-4}z\\x_{a+1}\ldots x_{t-1}y_1y_2\ldots y_{k+3}x_1$, a contradiction.  Assume that $a=t-1$.  From $t\geq 3$ and  $a=t-1$ it follows that  $
 A(\{x_1,x_2,\ldots , x_{t-2}\}\rightarrow 
 \{x_{t+1},x_{t+2},\ldots , x_{n-k-4}\})=\emptyset$.  This together with (13) implies that for some $s\in [1,t-2]$, $x_s\rightarrow x_t$. Since $d=k+3$ and $x_{j-1}\nrightarrow z$, from Claim 5(i) it follows that $z\rightarrow x_{s+1}$. Therefore, $C(z)=x_1x_2\ldots x_sx_{t}zx_{s+1}\ldots x_{t-1}x_{t+1}\ldots x_{n-k-4}y_1y_2\ldots y_{k+3}x_1$ is a Hamiltonian cycle in $D$, a contradiction.
 
Assume next that $t+1=n-d-1$. Recall that $b=t+1$ and $d^+(x_t,Y)=0$.
 
 Let $a\leq t-2$. If $x_{t-1}\rightarrow y_i$ with $i\in [1,2]$, then   $C(z)=x_1x_2\ldots x_ax_{n-d-1}zx_{a+1}\ldots x_{t-1}y_iy_2\ldots \\y_dx_1$ is a cycle of length at least $n-2$, a contradiction.  Therefore, we may assume that for every $i\in [1,2]$,  
 $d^-(y_i,\{x_{t-1}, x_t\})=0$.  This together with Claim 5(ii) implies that $d=k+4$, which in turn implies that (2) holds, in particular, $z\rightarrow x_t$.
  If $l=t-1$, then from $d^-(y_2,\{x_{t-1},x_t\})=0$ and Claim 5(ii) it follows that $y_2\rightarrow x_{a+1}$ and $C_{n-k-2}(z)=x_1x_2\ldots x_ax_{n-k-5}y_1y_2x_{a+1}\ldots  x_tzx_1$, a contradiction. Therefore,     we may assume that  $l\leq t-2$. It is easy to see that $a=t-2$
  (for otherwise $a\leq t-3$, $x_{t-2}\rightarrow y_1$ and  $C_{n-2}(z)=x_1x_2\ldots x_ax_{n-k-5}zx_{a+1}\ldots x_{t-2}y_1y_2\ldots y_{k+4}x_1$, a contradiction). 
  Using Claim 5(ii) and the facts that  $a=t-2\geq l\geq 2$, $d^-(y_1, \{x_{t-1},x_{t}\})=0$, 
 it is easy to see that $x_a\rightarrow y_1$.  
 From  $a=t-2\geq 2$ it follows that $d^-(x_{n-k-5}, \{x_{1},x_2,\ldots , x_{a-1}\})=0$. 
  Therefore by (13), there exist $s\in [1,a-1]$ and $j\in [t-1,t]$ such that $x_s\rightarrow x_j$.
   Then by (2), $z\rightarrow x_{s+1}$ and  $C(z)=x_1x_2\ldots x_sx_{j}x_tx_{n-k-5}zx_{s+1}\ldots x_a
   y_1y_2\ldots y_{k+4}x_1$ is a cycle of length at least $n-1$, a contradiction. 
Let now $a=t-1$. Recall that $b=t+1=n-d-1$. 
   Then from $a=t-1\geq 2$ we have that $d^-(x_{n-d-1}, \{x_{1},x_2,\ldots , x_{a-1}\})=0$.  This together with (13) implies  that for some $s\in [1,t-2]$, $x_s\rightarrow x_t$. 
   It is easy to see that 
   $z\nrightarrow x _{s+1}$ 
   (for otherwise, $z\rightarrow x_{s+1}$ 
   and  $C_n(z)=x_1x_2\ldots x_sx_{t}zx_{s+1}\ldots x_{t-1}x_{t+1} y_1y_2\ldots y_dx_1$,  a contradiction). From (2),  Claim 5(ii) and $z\nrightarrow x_{s+1}$ it follows that $d=k+3$ and $x_{t-1}\rightarrow y_1$. If $s\leq t-3$, then $z\rightarrow x_{s+2}$ and $C_{n-1}(z)=x_1x_2\ldots x_sx_{t}z x_{s+2}\ldots x_{t-1}x_{t+1} y_1y_2\ldots y_dx_1$,  a contradiction. Thus, we may assume that $s=t-2$. If $t-2\geq 2$, then we have that $A(\{x_1,x_2,\ldots , x_{t-2}\} \rightarrow \{x_t,x_{t+1}=x_{n-k-4}\})=\emptyset$. Therefore by (13), there is $p\in [1,t-3]$ such that
$x_p\rightarrow  x_{t-1}$ and $z\rightarrow  x_{p+1}$. 
If $l\leq t-2$, then $x_{t-2}\rightarrow  y_1$ 
and $C_{n}(z)=x_1x_2\ldots x_px_{t-1}x_t x_{t+1}zx_{p+1}\ldots x_{t-2} y_1y_2\ldots y_{k+3}x_1$, a contradiction. Assume that $l=t-1$. Then 
$y_1\rightarrow  x_{p+1}$ and $C_{n-k-2}=x_1x_2\ldots x_px_{t-1}x_{t+1}y_1x_{p+1}\ldots x_{t-2}x_tzx_1$, a contradiction.
Finally assume that $t-2=1$. Then $n-k-4=4$ and $d(x_3, Y)=0$. Therefore, $n+k\leq d(x_3)\leq 8$ and $n\leq 8$, which contradicts that $n\geq 9$.
 This completes the discussion of Subcase A.1.2. \\

\textbf{Subcase A.2.} $z\nrightarrow x_{a+1}$. 

From $z\nrightarrow x_{a+1}$, Claim 5(i), (1) and (2) it follows that $d=k+3$, $x_{b-1}\rightarrow z$ and
$$
\{x_{t},x_{t+1},\ldots , x_{n-k-4}\}\rightarrow z\rightarrow \{x_1,x_2,\ldots , x_a, x_{a+2},x_{a+3}, \ldots , x_t\}. \eqno (14)
 $$
  Assume first that
 $$
 A(\{x_1,x_2,\ldots , x_{t-2}\}\rightarrow \{x_{t+1},x_{t+2},\ldots , x_{n-k-4}\})=\emptyset.    \eqno (15)
 $$
 Then $a=t-1$, i.e., $x_{t-1}\rightarrow x_b$. Using (14), $d^+(z)\geq 2$ and Claim 2, we obtain that $t-1\geq 2$. From (13) and (15) it follows that there exists $s\in [1,t-2]$ such that $x_s\rightarrow x_t$. Then, since $z\rightarrow x_{s+1}$, $C_n(z)=x_1x_2\ldots x_{s}x_t\ldots x_{b-1}zx_{s+1}\ldots x_{t-1}x_b\ldots x_{n-k-4}y_1y_2\ldots y_{k+3}x_1$, a contradiction. 
 
 Assume next that (15) is not true. Then we may assume that
 $a\leq t-2$. Note that $z\rightarrow \{x_{a+2},\ldots , x_t\}$ (by (14)). If $y_i\rightarrow x_{a+1}$ with $i\in [1,k+3]$, then the cycle  $C(z)=x_1x_2\ldots x_{a}x_b\ldots x_{n-k-4}y_1\ldots y_ix_{a+1}\\ \ldots x_{b-1}zx_1$ has length at least $n-k-2$, a contradiction. We may therefore assume that $d^-(x_{a+1},Y)=0$. Let $b\geq t+2$. 
 Then, since $d=k+3$ and $t\geq l$, 
 from Claim 5(ii) 
 it follows that for some  $j\in [b-2,b-1]$, $x_j\rightarrow y_1$. 
 Then the cycle  $C(z)=x_1x_2\ldots x_{a}x_b\ldots x_{n-k-4}zx_{a+2}\ldots x_jy_1y_2\ldots y_{k+3}x_1$ has length at least $n-2$, a contradiction.
Let now $b=t+1$. We claim    that $l\leq t-1$. Assume that this is not the case, i.e., $l=t$. Then using Claim 5(ii) and the facts that $d=k+3$, $d^-(x_{a+1}, Y)=0$, we obtain that $x_t\rightarrow y_1$. Therefore, $C_{n-1}(z)=x_1x_2\ldots x_ax_{t+1}\ldots x_{n-k-4}zx_{a+2}\ldots x_{t}y_1\ldots y_{k+3}x_1$, a contradiction.
 This shows that $l\leq t-1$. From  $a\leq t-2$, $l\leq t-1$,
 $x_t\nrightarrow y_1$
and Claim 5(ii) it follows that $x_{t-1}\rightarrow y_1$. Therefore, if $a\leq t-3$, then $C_{n-2}(z) =x_1x_2\ldots x_ax_{t+1}\ldots x_{n-k-4}zx_{a+2}\ldots x_{t-1}y_1y_2\ldots y_{k+3}x_1$, a contradiction. We may therefore assume that $a=t-2$. 
 Assume first that $a\geq 2$. Since $a=t-2$, we have 
$$
 A(\{x_1,x_2,\ldots , x_{t-3}\}\rightarrow \{x_{t+1},x_{t+2},\ldots , x_{n-k-4}\})=\emptyset.    
 $$
This together with (13) implies that there exist $s\in [1,a-1=t-3]$ and $p\in [t-1,t]$ such that $x_s\rightarrow x_p$. Then by (14), the cycle $C(z)=x_1x_2\ldots x_sx_px_tzx_{s+1}\ldots x_{t-2}x_{t+1}\ldots x_{n-k-4}y_1y_2\ldots  y_{k+3}x_1$ has length at least $n-1$, a contradiction. 

 Assume next that  $a=1$. Then $t=3$. Let $t+1\leq n-k-5$. Since $b=t+1$, we have $d^+(x_{1}, \{x_{t+2},x_{t+3},\ldots ,  x_{n-k-4}\})=0$. Again using (13), we obtain that  there exist $p\in [t-1,t]$ and $q\in [t+2,n-k-4]$ such that 
$x_p\rightarrow x_q$. Recall that $z\rightarrow x_t$ and $x_{q-1}\rightarrow \{z,y_1\}$.
Therefore, if $p=t$, then $C_{n-1}(z)=x_1x_{t+1}\ldots x_{q-1}zx_tx_q\ldots 
x_{n-k-4}y_1y_2\ldots y_{k+3}x_1$, and if $p=t-1$, then  
$C_{n}(z)=x_1x_2\ldots x_{t-1}x_q\ldots x_{n-k-4}zx_t\ldots x_{q-1}y_1y_2\ldots y_{k+3}x_1$. Thus, in both cases, we have a contradiction.
This completes the discussion of Subcase A.2, and also completes the proof of the theorem when  $D\langle Y\rangle$ is  strong.\\

\textbf{Case B.} $D\langle Y\rangle$ is not strong.

Since $y_1y_2\ldots y_{d}$ is a path in $D\langle Y\rangle$ and $k+3\leq d\leq k+4$, using the fact that every vertex $y_i$
with $i\in [1,d]$ cannot be inserted into $P$ (Claim 4) and Lemma 3.2, we obtain $d(y_i, V(P))\leq n-d$, $d(y_i, Y)\geq d+k$. Now, we claim that $d=k+4$ and $k=0$. Indeed, since $D\langle Y\rangle$ is not strong, $y_1y_2\ldots y_d$ is a Hamiltonian path in $D\langle Y\rangle$ and $d(y_i)\geq n+k$, it follows that for some $l\in [2,d-1]$, $y_l\rightarrow y_1$ and $d^-(y_1,\{y_{l+1},y_{l+2},\ldots , y_d\})=d^+(y_d,\{y_{1},y_2\ldots , y_l\})=0$. From this we have $k+d\leq d(y_1,Y)\leq d-l+2(l-1)=d+l-2$ and $k+d\leq d(y_d,Y)\leq l+2(d-l-1)=2d-l-2$. Therefore, $k\leq l-2$ and $d\geq k+l+2$. From the last two inequalities and the facts that $d\leq k+4$, $l\geq 2$ it follows that $d=k+4$ and $k=0$. 
Therefore, $d(y_i, V(P))\leq n-4$ and $d(y_i, Y)\geq 4$.
 Since $D\langle Y\rangle$ is not strong and $y_1y_2y_3y_4$ is a path in $D\langle Y\rangle$, it is not difficult to check that for all $i\in [1,4]$, $d(y_i,Y)=4$, $d(y_i,V(P))=n-4$, the arcs $y_1y_3, y_1y_4,y_2y_1,  y_2y_4,y_4y_3$ also are in $A(D)$ and $A(\{y_3,y_4\}\rightarrow \{y_1,y_2\})=\emptyset$. Since $D$ has no $C(z)$-cycle of length at least $n-2$ and any vertex $y_i$ with $i\in [1,4]$ cannot be inserted into $P=x_1x_2\ldots x_{n-5}$, using Lemma 3.3 and Proposition 3, it is not difficult to show that  there are two integers $l_1$ and $l_2$ with $2\leq l_1, l_2\leq n-6$ such that
$$
 \left\{ \begin{array}{lc}\{x_{l_1},\ldots ,x_{n-5}\}\rightarrow \{y_1,y_2\}\rightarrow \{x_1,\ldots , x_{l_1}\}, \\  \{x_{l_2},\ldots ,x_{n-5}\}\rightarrow \{y_3,y_4\} \rightarrow \{x_1,\ldots , x_{l_2}\}.   \\ \end{array} \right.  \eqno (16)
$$
It is easy to see that $l_1\geq l_2$. Indeed, if $l_1\leq l_2-1$, then from (16) it follows that  $x_{l_1}\rightarrow y_1$, $y_4\rightarrow x_{l_1+1}$ and hence, $C_n(z)=x_1\ldots x_{l_1}y_1y_2y_3y_4x_{l_1+1}\ldots x_{n-5}zx_1$, a contradiction.   Since $D$ is 2-strong,  (16) together with (2) implies that there are two integers $p\in [1,l_2-1]$ and $q\in [l_2+1,n-5]$ such that $x_p\rightarrow x_q$ (for otherwise $D-x_{l_2}$ is not strong). Assume first that $l_2\leq t$. 
Then from (2) and (16), respectively, we have $z\rightarrow x_{p+1}$ and  $x_{q-1}\rightarrow y_3$. Therefore,  $C_{n-2}(z)=x_1\ldots x_px_q\ldots x_{n-5}zx_{p+1}\ldots x_{q-1}y_3y_4x_1$, a contradiction.
Assume next that $l_2\geq t+1$. Then by (16),   $y_4\rightarrow x_{p+1}$,  and by (2), $x_{q-1}\rightarrow z $. Therefore,
$C(z)=x_1x_2\ldots x_px_q\ldots x_{n-5}y_1\ldots y_4x_{p+1}\ldots x_{q-1}zx_1$ is a Hamiltonian cycle in $D$, which is contradiction. This completes the discussion of Case B. 
The theorem is proved. \end{proof}

\acknowledgements

I am grateful to Professor Gregory Gutin for motivating me to present the complete proof of Theorem 1.7, and very  thankful to the anonymous referees for a careful reading of our manuscript and for their detailed and helpful comments that have improved the presentation. Also thanks to Dr. Parandzem Hakobyan for formatting the manuscript of this paper.


\begin{thebibliography}{30}
\providecommand{\natexlab}[1]{#1}
\providecommand{\url}[1]{\texttt{#1}}
\expandafter\ifx\csname urlstyle\endcsname\relax
  \providecommand{\doi}[1]{doi: #1}\else
  \providecommand{\doi}{doi: \begingroup \urlstyle{rm}\Url}\fi

  

\bibitem [Adamus(2017)]{[1]} J.~Adamus. \newblock A  degree sum condition for hamiltonicity in balanced bipartite digraphs. \newblock \emph{ Graphs and Combinatorica}, 33(1):\penalty0 43--51, 2017.  doi: /10.1007/s00373-016-1751-6.

\bibitem [Adamus(2021)]{[2]} J.~Adamus.  \newblock On dominating pair degree conditions for hamiltonicity in balanced bipartite digraphs.\newblock \emph { Discrete Math.}, 344(3):\penalty0 112240, 2021.  doi: 10.1016/j.disc.2020.112240. 



\bibitem [ Adamus and  Adamus (2012)]{[3]} J.~Adamus and L.~Adamus. \newblock  A degree codition for cycles of maximum length in bipartite digraphs. \newblock \emph {Discrete Math.}, 312(6):\penalty0 1117--1122, 2012.   doi: 10.1016/j.disc.2011.11.032.

\bibitem [Adamus, Adamus and Yeo (2014)]{[4]} J.~Adamus, L.~Adamus and A.~Yeo. \newblock On the Meyniel condition for hamiltonicity  in bipartite digraphs. \newblock \emph { Discrete Math. and Theoretical Computer Science}, 16(1):\penalty0 293-302, 2014. doi:  10.46298/dmtcs.1264.

\bibitem[Bang-Jensen and Gutin (Springer-Verlag, London, 2000)]{[5]}
J.~Bang-Jensen and G.~Gutin.
\newblock \emph{Digraphs: Theory, Algorithms and Applications},
  Springer-Verlag, London, 2000.

\bibitem[Bermond and Thomassen(1981)]{[6]}
J.-C. Bermond and C.~Thomassen. 
\newblock Cycles in digraphs $-$ a survey.
\newblock \emph{J. Graph Theory}, 5\penalty0 (1):\penalty0 1--43, 1981.  doi: 10.1002/jgt.3190050102.

\bibitem[Bondy and Thomassen(1977)]{[7]}
J.A.~Bondy and C.~Thomassen.

\newblock A short proof of Meyniel's theorem.
\newblock \emph{Discrete Math.}, 19:\penalty0 195--197, 1977.
doi:  10.1016/0012-365X(77)90034-6.

\bibitem[Darbinyan(1982)]{[8]} S.Kh. ~Darbinyan. \newblock Cycles of any length in digraphs with large semidegrees. \newblock \emph{ Aakd. Nauk  Armyan. SSR Dokl.}, 75\penalty0(4):\penalty0 147--152, 1982.   

 \bibitem[Darbinyan(1986)]{[9]} S.Kh.~Darbinyan. \newblock A sufficient condition for the Hamiltonian property of  digraphs   with large semidegrees. \newblock \emph{Aakd. Nauk  Armyan. SSR Dokl.}, 82\penalty0 (1):\penalty0 6--8, 1986.  
   
  \bibitem[Darbinyan(1990)]{[10]}
S.Kh.~Darbinyan.
 \newblock Hamiltonian and strongly Hamilton-connected digraphs.
\newblock \emph{Akad. Nauk Armyan. SSR Dokl.}, 91\penalty0 (1):\penalty0 3--6, 1990 (for a detailed proof see arXiv:
  1801.05166v1, 16 Jan. 2018).
   
  \bibitem[Darbinyan(1990a)]{[11]} S.Kh.~Darbinyan. \newblock A sufficient condition for a digraph to be Hamiltonian. \newblock \emph {Aakd. Nauk  Armyan. SSR Dokl.}, 91(2):\penalty0  57--59, 1990a. 
  
 \bibitem[Darbinyan(2022)] {[12]} S.Kh.~Darbinyan. \penalty0 On an extension of the Ghouila-Houri theorem,  \newblock \emph{Mathematical Problems of Computer Science}, 58:\penalty0 20--31, 2022. doi: 10.51408/1963-0089.

 \bibitem[Darbinyan(2024)] {[13]} S.Kh. ~Darbinyan. \newblock On Hamiltonian cycles in a 2-strong digraphs with large degrees and cycles,  \newblock \emph{Patern Recognition and Image Analaysis}, 34:\penalty0(1) 62--73, 2024. doi: 10.1134/S105466182401005X .

\bibitem[Ghouila-Houri(1960)]{[14]}
A.~Ghouila-Houri. 
\newblock Une condition suffisante d'existence d'un circuit hamiltonien.
\newblock \emph{C. R. Acad. Sci. Paris Ser. A-B}, 251:\penalty0 495--497, 1960.

\bibitem [Goldberg,  Levitskaya and  Satanovskyi(1971)] {[15]}   M.K.~Goldberg, L.P. ~Levitskaya and L.M. ~Satanovskyi. \newblock On one strengthening of the Ghouila-Houri theorem. \newblock \emph { Vichislitelnaya Matematika i Vichislitelnaya Teknika}, 2:\penalty0 56-61, 1971.

\bibitem [Gould(2014)]{[16]} R.~Gould. \newblock Resent Advances on the Hamiltonian Problem: Survey III. \newblock \emph{ Graphs and Combinatorics}, 30:1-46, 2014. 
doi:  10.1007/s00373-013-1377-x.

\bibitem[H\"{a}ggkvist and Thomassen(1976)]{[17]}
R.~H\"{a}ggkvist and C.~Thomassen.
 \newblock On pancyclic digraphs.
\newblock \emph{J. Combin. Theory Ser. B}, 20:\penalty0 20--40, 1976.  doi:  10.1016/0095-8956(76)90063-0.

\bibitem[K\"{u}hn and Osthus(2012)]{[18]}
D.~K\"{u}hn and D.~Osthus.
 \newblock A survey on Hamilton cycles in directed graphs.
\newblock \emph{European J. Combin.}, 33(5):\penalty0 750--766, 2012.  doi:  10.1016/j,ejc.2011.09.030.

\bibitem[Meyniel(1973)]{[19]}
M.~Meyniel.
\newblock Une condition suffisante d'existence d'un circuit hamiltonien dans un   grape oriente.
\newblock \emph{J. Combin. Theory Ser. B}, 14:\penalty0 137--147, 1973.   doi:  10.1016/0012-365X(77)90034-6.

\bibitem[Nash-Williams(1969)]{[20]}
C.~St. J.~A. Nash-Williams. \newblock
Hamilton circuits in graphs and digraphs, the many facets of graph
  theory.
\newblock \emph{Springer Verlag Lecture Notes}, 110:\penalty0 237--243, 1969. 

\bibitem[ Overbeck-Larisch(1976)] {[21]} M. ~Overbeck-Larisch. \newblock Hamiltonian paths in oriented graphs. \newblock \emph {J. Combin. Theory Ser. B}, 21:\penalty0 76--80, 1976. doi: 10.1016/0095-8956(76)90030-7.

\bibitem [Thomassen (1981)]{[22]} C.~Thomassen.  \newblock  Long cycles in digraphs,  \newblock \emph{ Proc. London Math. Soc.}, (3)42: \penalty0 231--251, 1981. doi:  10.1112/PLMS/S3-42.2231.

\bibitem [Wang (2021)]{[23]} R.~Wang. \newblock Extremal digraphs on Woodall-type condition for hamiltonian cycles  in balanced bipartite digraphs, \newblock \emph{J. Graph Theory},  97(2):\penalty0 194--207, 2021.  doi: 10.1002/jgt.22649.


\bibitem [Wang (2022)] {[24]} R.~Wang. \newblock A note on dominating pair degree condition for hamiltonian cycles  in balanced bipartite digraphs, \newblock \emph{ Graphs and Combinatorics},  38: \penalty0 13, 2022. doi: 10.1007/s00373-021-02404-8. 

\bibitem [Wang and  Wu(2021)] {[25]} R.~Wang and L.~Wu. \newblock Cycles of many lengths  in balanced bipartite digraphs on dominating and dominated degree condition, \newblock \emph{Discussiones  Mathematicae Graph Theory}, \penalty0 1--23, 2021. doi: 10.7151/dmgt.2442.

\bibitem [Wang,  Wu and  Meng(2022)]{[26]}   R.~Wang, L.~Wu and W. ~Meng. \newblock Extremal digraphs on Meyniel-type condition for hamiltonian cycles  in balanced bipartite digraphs. \newblock \emph {Discrete Math. and Theoretical Computer Science}, 23:\penalty0 (3),  2022.
doi: 10.46298/dmtcs.5851.
 

\bibitem[Woodall(1972)]{[27]}
D.~Woodall.
\newblock Sufficient conditions for circuits in graphs.
\newblock \emph{Proc. London Math. Soc.}, 24:\penalty0 739--755, 1972.  doi: 10.1112/plms/s3-24.4.739.

\bibitem [Zhang,  Zhang and Wen(2013)] {[28]}  Z.-B.~Zhang, X. ~Zhang and X. ~Wen. \newblock Directed Hamilton cycles in digraphs and matching alternating Hamilton cycles in bipartite graphs. \newblock \emph {SIAM J. on Discrete Math.}, 27(1):\penalty0 274--289, 2013. doi: 10.1137/110837188.



\end{thebibliography}
\end{document}